\documentclass[10.5pt]{article}

\usepackage{amsmath}
\usepackage{amscd}
\usepackage{amssymb}
\usepackage{rotating}
\usepackage{amsfonts}
\usepackage{amsthm}
\usepackage{verbatim}
\usepackage{stmaryrd}
\usepackage[top=1in, bottom=1in, left=1.25in, right=1.25in]{geometry}

\usepackage{titlesec}
       \titleformat{\chapter}[display]
             {\normalfont\Large\bfseries}{\thechapter}{11pt}{\Large}
       \titleformat{\section}
             {\normalfont\large\bfseries}{\thesection}{10.5pt}{\large}
       \titlespacing*{\chapter}{0pt}{0pt}{15pt} 
       \titlespacing*{\section}{0pt}{3.5ex plus 1ex minus .2ex}{2.3ex plus .2ex}

    \usepackage[titletoc]{appendix}

\input xy
\xyoption{all}

\allowdisplaybreaks[1]
\newcommand{\pd}{\partial}

\newcommand{\cM}{{\mathcal M}}

\newcommand{\lbr}{{\llbracket}}
\newcommand{\rbr}{{\rrbracket}}

\DeclareMathOperator{\ch}{ch}
 
 \DeclareMathOperator{\ev}{ev}

\newtheorem{theorem}{Theorem}[section]
\newtheorem{theorem/definition}{Theorem/Definition}[section]

\newtheorem{proposition}{Proposition}[section]
\newtheorem{lemma}{Lemma}[section]

\newtheorem{corollary}{Corollary}[section]
\newtheorem{Conjecture}{Conjecture}

\theoremstyle{remark}
\newtheorem{remark}{Remark}[section]

\theoremstyle{definition}

\newcommand{\be}{\begin{equation}}
\newcommand{\ee}{\end{equation}}
\newcommand{\bea}{\begin{eqnarray}}
\newcommand{\eea}{\end{eqnarray}}
\newcommand{\ben}{\begin{eqnarray*}}
\newcommand{\een}{\end{eqnarray*}}
\newcommand{\bet}{\begin{equation}
\begin{split}}
\newcommand{\eet}{\end{split}
\end{equation}}

\begin{document}
\title
{\Large On the Crepant Resolution Conjecture for Gromov-Witten Gravitational Ancestors in All Genera for Surface Singularities
\thanks{Date: December 14, 2013}}
\author{\normalsize Xiaowen Hu}
\date{}
\maketitle

\begin{abstract}
We state a version of the crepant resolution conjecture for total ancestor potentials for surface singularities, and reduce the conjecture to the quantum McKay correspondence conjecture of J.Bryan and A.Gholampour and a vanishing conjecture for hurwitz-hodge integrals. In particular, for singularities of type A, we prove the conjecture. We also suggest an approach towards a proof for the general cases by Teleman's reconstruction theorem for semi-simple cohomology field theories.
\end{abstract}

\section{Introduction}
Let $\mathcal{X}$ be an effective orbifold with the coarse moduli space $X$ and $Y\rightarrow X$ be a crepant resolution. The general principle of McKay correspondence expects that the geometry of $\mathcal{X}$ coincides with that of $Y$. This principle is extended to a quantum version, and has been stated as the crepant resolution conjecture for Gromov-Witten invariants. For the history of this conjecture, we refer the reader to \cite{CoatesRuan}, \cite{Zhou0811} and more references there. Here we start from a conjecture of J.Bryan and A.Gholampour \cite{BGh0707}.

\begin{Conjecture}\label{13}
Let  $F_{0}^{ \widehat{\mathbb{C}^{2}/G}}(y_{0},y_1,\cdots,y_n;\mathbf{q})$, $F_{0}^{ [\mathbb{C}^{2}/G]}(x_{0},x_1,\cdots,x_n)$  denote  the  $\mathbb{C}^{*}$-equivariant  genus 0 orbifold Gromov-Witten potential  of $\widehat{\mathbb{C}^{2}/G}$ and the orbifold $\mathcal{X}= [\mathbb{C}^{2}/G]$. Then after the change of variables
\bea\label{12}
y_{0}&=&x_{0},\nonumber\\
y_{R}&=&\frac{1}{|G|}\sum_{g\in G}\sqrt{2-\chi_{\rho_{1}}(g)}\overline{\chi}_{R}(g)x_{\lbr g\rbr},\\
q_{R}&=&\exp\Big(\frac{2\pi\sqrt{-1}\dim R}{|G|}\Big),\nonumber
\eea
we have $F_{0}^{ \widehat{\mathbb{C}^{2}/G}}=F_{0}^{ [\mathbb{C}^{2}/G]}$.
\end{Conjecture}
Note that as a formal series of the $y_{i}$'s, the coefficients of  $F_{0}^{ \widehat{\mathbb{C}^{2}/G}}(y_{0},y_1,\cdots,y_n;\mathbf{q})$ is not convergent at $q_{R}=\exp\Big(\frac{2\pi\sqrt{-1}\dim R}{|G|}\Big)$. Thus we need to make precise what we mean by this change of variables of $q_{R}$ (e.g., if we want to make precise predictions for hurwitz-hodge integrals via the crepant resolution conjecture). In fact in the derivation of the prediction for $F_{0}^{ [\mathbb{C}^{2}/G]}$ by the above conjecture in \cite{BGh0707}, the following convention was implicit.\\
\textbf{Convention}: The coefficients of $F_{0}^{ \widehat{\mathbb{C}^{2}/G}}(y_{0},y_1,\cdots,y_n;\mathbf{q})$, viewed as functions with complex variables $q_{R}$'s, are well-defined near the origin, and we analytically continuate them along the ray from $0$ to $\exp\Big(\frac{2\pi\sqrt{-1}\dim R}{|G|}\Big)$ and take values over there.\\

 In this article, we shall always use this convention. \\

This conjecture is proved in \cite{BGr} for $G=\mathbb{Z}_{2}$, by a result of \cite{FP}. In \cite{CCIT2}, a version of the above conjecture
for type A groups was proved, as a corollary of the corresponding correspondence of the $J$-functions. However, in that correspondence the variables $q_{R}$ and $y_{R}$ are mixed in a way such that it seems (at least to the author) one cannot easily deduce the conjecture in the precise form above. Fortunately, later this conjecture for the type A groups is included in the main theorem of \cite{Zhou0811}; in fact, its main theorem is  stated for the correspondence of  stationary \emph{reduced} Gromov-Witten invariants (combined with the work of  \cite{Maulik}) for any higher genera,
 and when restricted to genus zero, the proof would be greatly simplified. Based on the validity of the formula of $F_{0}^{ [\mathbb{C}^{2}/G]}$  for groups of type A, the author proved (see \cite{Hu}) the conjecture in \cite{BGh0707} on $F_{0}^{ [\mathbb{C}^{2}/G]}$ for groups of type D, thus established the above conjecture in this case.\\

 It is natural to extend the above correspondence to higher genera and gravitational ancestor or descendant invariants, not only for the reduced theory ($=$ the $t$-linear part of the full equivariant theory). It is reasonable to expect that, the Gromov-Witten correlators on the \emph{orbifold-side} can be obtained by those on the \emph{resolution-side} through \emph{summing the contracted curve classes}. Y.P. Lee et al. have made some observations in \cite{LLW}. In this article, we consider the surface singularities. On the resolution side, the primary invariants seems less interesting; by
  the divisor equation and the vanishing result on the top components of the equivariant virtual cycles \cite{OP}, these correlators reduce to the degree zero ones and the classical invariants, as in \cite{BGh0707}. On the orbifold side, by dimensional reason, the difficulty is concentrated in genus one, where the conjectural correspondence gives a prediction on the vanishing of Hurwitz-Hodge integrals, which is similar to the vanishing result mentioned above. We state this prediction as a conjecture (see conjecture \ref{7}) , and give a proof for the cases of only one marked point.\\

  For ancestor invariants, a more general vanishing conjecture (see conjecture \ref{9}) is stated. Assuming the validity of this conjecture and the conjecture \ref{13}, the argument of \cite{FSZ} applies, and therefore implies our main theorem.
\begin{theorem}\label{10}
Assuming the quantum McKay correspondence conjecture \ref{13}, and the conjecture \ref{9}, the total ancestor Gromov-Witten potential $\mathcal{A}^{\widehat{\mathbb{C}^{2}/G}}$, after the analytic continuation and the change of variables (\ref{12}), equals the total ancestor Gromov-Witten potential $\mathcal{A}^{[\mathbb{C}^{2}/G]}$.
\end{theorem}
 Since the conjecture \ref{9} for groups of type A follows plainly from an orbifold version of Mumford's relations (prop.3.2 of \cite{BGP}), we have \begin{corollary}\label{14}
 The total ancestor Gromov-Witten potential $\mathcal{A}^{\widehat{\mathbb{C}^{2}/\mathbb{Z}_{n}}}$, after the analytic continuation and the change of variables (\ref{13}), equals the total ancestor Gromov-Witten potential $\mathcal{A}^{[\mathbb{C}^{2}/\mathbb{Z}_{n}]}$.
\end{corollary}
As a biproduct, we recover a result in \cite{Zhou0811} for higher genera, see remark \ref{14}. \\

While the proofs in this article are routine computations case by case, the author would like to make some further comments.\\

(i) As A.Givental conjectured in \cite{Givental1} and \cite{Givental2}, and C.Teleman proved in \cite{Teleman}, a \emph{semi-simple homogeneous} cohomological field theory (and thus the associated total ancestor potential) is determined by the Frobenius manifold defined by the genus zero primary potential. The equivariant Gromov-Witten potential considered in our cases are (generically) semi-simple but not homogeneous, and A.Givental and C.Teleman's results showed that in these cases the total ancestor potential is determined by the genus zero primary potential \emph{up to an ambiguity} (called \emph{Hodge twisting} in \cite{Teleman}). This ambiguity can be viewed as a series of functions depending on the equivariant parameters $t$ and the Novikov variables $\mathbf{q}$, but not on the cohomological variables. Thus it is very natural to expect that this ambiguity is determined by the additional data  $F_{1,1}$ and $F_{g,0}$, for all $g\geq 2$. In particular, we believe that the vanishing conjecture \ref{9} can be proved in this way, and this is why we verify  the conjecture \ref{7} for $F_{1,1}$. \\

(ii) As we mentioned earlier, the result on the correspondence of (descendant) $J$-functions in \cite{CCIT2} seems not easy to deduce the conjecture \ref{13} (for groups of type A) in the precise form. In Appendix B, we  show that for $\widehat{\mathbb{C}^{2}/\mathbb{Z}_{2}}$,
the \emph{small} quantum product and the \emph{ancestor divisor equation}  will imply that  the ancestor $J$-function will satisfy a hypergeometric differential equation after change of variables and analytic continuations. Together with the result on $I$-function of
$[\mathbb{C}^{2}/\mathbb{Z}_{2}]$ (\cite{CCIT2}), we prove the correspondence for ancestor $J$-functions and thus the conjecture \ref{13} in this case. The result is of course not new\footnote{From the  correspondence for primary invariants, and the theorem 4.3 in Chap.III of \cite{Manin}, one can deduce the correspondence for ancestor invariants in genus zero.}, but may have its own interest, and we hope that this approach can be extended to other circumstances, e.g., the \emph{non-hard Lefschetz} crepant resolution conjecture, the LG/GW correspondence...  \\

(iii) Until now we have been putting the \emph{descendant} invariants aside. Thus now a natural question is, what is the relation between the correspondence of ancestor $J$-functions (resp., the total ancestor potential) and that of the descendant $J$-functions (resp. the total descendant potential)? The author will not investigate this aspect in this article, and refer the reader to \cite{Givental1}, \cite{Givental2}, and \cite{CoatesRuan}.\\

\textit{Acknowledgements.} The author thanks Prof.~Jian Zhou for his great patience and guidance during all the time. He also thanks Yutao Ding, Di Yang, Xiaobo Zhuang, Zhilan Wang, and Hanxiong Zhang for helpful discussions, and thanks Prof.  Hsian-Hua Tseng for helpful communications.

\section{Preliminaries on the equivariant Gromov-Witten theory}
In this section we recall some rudiments on equivariant Gromov-Witten theory. For more details, see \cite{BF}, \cite{Beh}, \cite{GP}, and
\cite{AGV1}, \cite{Ts}.
Let $T$ be a complex torus, $V$ be a nonsingular complex quasi-projective variety with a $T$-action. Let $\overline{\cM}_{g,n}(V, \beta)$ be the moduli space of stable maps. By the existence of a $T$-equivariant polarization (\cite{MFK}, theorem 1.7, chap.1), we have a $T$-equivariant perfect obstruction theory on  $\overline{\cM}_{g,n}(V, \beta)$, thus obtain a $T$-equivariant virtual fundamental cycle in the equivariant
 chow group $A_{\star}^{T}(\overline{\cM}_{g,n}(V, \beta))$. When $V$ is projective, we have the the corresponding \emph{system of equivariant Gromov-Witten invariants}
\ben
I_{\tau}(\beta): H^{\star}_{T}(V)^{\otimes S_{\tau}}\rightarrow H^{\star}(\overline{\cM}_{\tau})\otimes R_{T},
\een
for every stable modular graph $\tau$ and every effective cycle marking $\beta$ of $\tau$. Here $R_{T}=\text{Sym}_{\mathbb{Q}}(\hat{T})$, $\hat{T}$ is the character group of $T$.
\begin{theorem}\label{21}
The system of equivariant Gromov-Witten invariants a nonsingular projective variety $V$ satisfy the usual axioms of usual Gromov-Witten invariants, and also the \emph{tautological relations} induced by the tautological relations on the moduli spaces of stable curves.
\end{theorem}
Proof: The proof goes  the same way as in \cite{Beh}. Only the divisor axiom needs some words. To define the pairing of a equivariant divisor class and the (usual) effective cycle class $\beta$, we need the fact that one can choose a $T$-stable cycle representing $\beta$, and the pairing is independent of the choice of the cycle. This fact follows from theorem 1 in \cite{FMSS} (see also \cite{Hirshowitz}).\hfill\qedsymbol\\

\begin{remark}
$R_{T}$ is viewed as a graded ring, with every character of (real) degree two, such that the homogeneous axiom holds. One can also let the characters take varying values in $\mathbb{C}$, thus obtain an \emph{inhomogeneous cohomological field theory} and the reconstruction result (up to an ambiguity) of \cite{Teleman} applies when the Frobenius manifold defined by the genus zero equivariant Gromov-Witten invariants  is generically semisimple.\\
\end{remark}
Sometimes we need to allow $V$ to be quasi-projective, with the following assumption :\\
Assumption ($\clubsuit$): There exists a closed subvariety $W$ of $V$,  such that every nonzero effective cycle $\beta$ lies in $W$, and $W$ is projective itself. Moreover, the $T$-fixed locus of $V$ lies in $W$.\\

In the following in this article, when we talk about the equivariant Gromov-Witten invariants of  a nonsingular quasi-projective variety $V$, we always assume the assumption ($\clubsuit$). Thus when $\beta$ is nonzero, $\overline{\cM}_{g,n}(V, \beta)$ is still proper, and we can push forward the virtual fundamental class. When $\beta=0$, $\overline{\cM}_{g,n}(V, \beta)\cong\overline{\cM}_{g,n}\times V$ is non-proper, but the $T$-fixed locus lies in $\overline{\cM}_{g,n}\times W$, so we can push forward the \emph{localized} virtual fundamental class. In this way, the \emph{system of local Gromov-Witten invariants} takes values in $H^{\star}(\overline{\cM}_{\tau})\otimes \mathcal{Q}_{T}$, where $\mathcal{Q}_{T}$ denotes the ring $R_{T}$ localized at all homogeneous elements. In summary, we have
\begin{theorem}
The system of local equivariant Gromov-Witten invariants of a nonsingular quasi-projective variety $V$ satisfy the usual axioms of usual Gromov-Witten invariants, and also the \emph{tautological relations} induced by the tautological relations on the moduli spaces of stable curves.
\end{theorem}\hfill\qedsymbol\\

Now we recall some standard notations. Let $\gamma_{1},\cdots,\gamma_{m}$ be a basis of $H_{T}^{*}(V)$, with the corresponding variables $x_{1},\cdots,x_{m}$, then the genus $g$, $n$-point primary Gromov-Witten potential is defined to be
\ben
F_{g,n}^{V}=\sum_{d}\frac{1}{n!}\langle \gamma^{n}\rangle_{g,n,d}^{V}q_{1}^{d_{1}}\cdots q_{r}^{d_{r}},
\een
where $\gamma=\gamma_{1}x_{1}+\cdots+\gamma_{m}x_{m}$, and $d=d_{1}\beta_{1}+\cdots+d_{r}\beta_{r}$ runs over the effective classes in $H_{2}(V)$. Note that the correlators not in the stable range is defined to be zero by convention. The genus $g$ primary Gromov-Witten potential is defined to be
\ben
F_{g}^{V}=\sum_{n\geq 0}F_{g,n}^{V}.
\een
We denote by $\psi_{i}$ the first chern class of the universal cotangent line bundle over $\overline{\cM}_{g,n}(V, \beta)$ at the $i$-th marked point, and by $\phi_{i}$ the pull backed  first chern class of the universal cotangent line bundle over $\overline{\cM}_{g,n}$ by the stabilization morphism when $2g-3+n\geq 0$, i.e., in the \emph{absolute} stable range. We use curved symbol $\mathcal{F}_{g,n}^{V}$, $\mathcal{F}_{g}^{V}$ to denote the ancestor potentials, defined in a similar way as above, and the \emph{total ancestor potential} is defined to be
\ben
\mathcal{A}^{V}=\exp \sum_{g\geq 0}\hbar^{g-1}\mathcal{F}_{g}^{V}.
\een

\begin{remark}
The above definitions and properties remain valid for a $T$-equivariant orbifold $\mathcal{X}$ with quasi-projective coarse moduli. $H_{T}^{*}(V)$ should be replaced by the Chen-Ruan cohomology group $H_{T}^{*}(I\mathcal{X})$, where $I\mathcal{X}$ is the inertia stack of $\mathcal{X}$, while the divisor axiom still holds only for the usual $H_{T}^{2}(\mathcal{X})$. The moduli stacks of stable maps are replaced by $\overline{\cM}_{g,n}(\mathcal{X}, \beta)$ (see \cite{AGV1}, \cite{Ts}). The other necessary modifications are obvious.
\end{remark}

\section{Local equivariant primary Gromov-Witten invariants of ADE-resolutions in higher genera}
For a finite subgroup $G\subset SL(2,\mathbb{C})$, let $\widehat{\mathbb{C}^{2}/G}$ be the (unique) crepant resolution of $\mathbb{C}^{2}/G$. There is a natural torus action of $T$ on $\widehat{\mathbb{C}^{2}/G}$ induced by the Hilbert scheme description of $\widehat{\mathbb{C}^{2}/G}$: when $G$ is of type A, $T= \mathbb{C}^{*}\times  \mathbb{C}^{*}$; when $G$ is of type D or type E, $T= \mathbb{C}^{*}$. It is obvious that $\widehat{\mathbb{C}^{2}/G}$ satisfies the assumption ($\clubsuit$). Consider the moduli space of stable maps $\overline{\cM}_{g,n}(\widehat{\mathbb{C}^{2}/G},d)$. When $d>0$, it is compact, and  the components in top (cohomological) degree of the equivariant virtual fundamental cycle agree with the virtual cycle in the non-equivariant Gromov-Witten theory, thus has dimension $g-1+n$. However, as shown in \cite{OP}, based on the results of \cite{Ran} and \cite{Manetti}, the top degree components are trivial. In the primary Gromov-Witten invariants for these spaces, every marked point can carry at most degree one cohomological classes, therefore by dimensional reason the primary invariants for $d>0$ and $g>0$ are all zero. Since the genus zero theory has been treated in \cite{BGh0707}, we only need to consider the degree zero cases.\\

For general target spaces $V$, We have $\overline{\cM}_{g,n}(v,0)\cong \overline{\cM}_{g,n}\times V$, and
\ben
[\overline{\cM}_{g,n}(V,0)]_{T}^{\text{vir}}=e_{T}(T_{V}\boxtimes \mathbb{E}^{\vee})\cap ( [ \overline{\cM}_{g,n}]\times [V]),
\een
where $e_{T}$ denotes the equivariant euler class, and $\mathbb{E}^{\vee}$ the dual Hodge bundle over $ \overline{\cM}_{g,n}$. To define the equivariant Gromov-Witten invariants, we need to localize this cycle to the fixed loci of the $T$-action.

\begin{theorem}\label{4}
The only nonzero equivariant primary gromov-witten invariants for $Y_{G}$ with degree zero and genus $> 0$ are
\bea\label{5}
\langle 1\rangle_{1,1,0}^{Y_{G}}=\left \{\begin{array}{ll}
-\frac{1}{24|G|}\frac{t_{1}+t_{2}}{t_{1}t_{2}}, & \text{if}\hspace{0.16cm} G\hspace{0.14cm}\text{is of type A};\\
-\frac{1}{12|G|t},& \text{if}\hspace{0.16cm} G\hspace{0.14cm}\text{is of type D or E},\end{array}\right.
\eea
and
\bea\label{6}
\langle \cdot\rangle_{2,0,0}^{Y_{G}}=\left \{\begin{array}{ll}
-\frac{1}{5760|G|}\frac{t_{1}+t_{2}}{t_{1}t_{2}}, & \text{if}\hspace{0.16cm} G\hspace{0.14cm}\text{is of type A};\\
-\frac{1}{2880|G|t},& \text{if}\hspace{0.16cm} G\hspace{0.14cm}\text{is of type D or E}.\end{array}\right.
\eea
\end{theorem}

\emph{Proof}: We first give a detailed computation for $G=\hat{E}_{6}$. We use the following graph to indicate the fixed loci of $Y_{\hat{E}_{6}}$,
$$ 
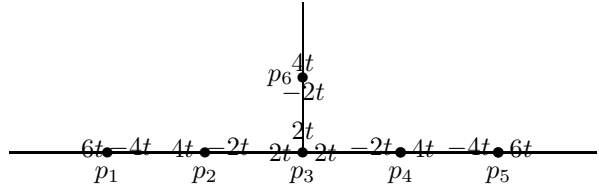
\xy
(0,0); (13,0), **@{-};(26,0), **@{-}; (39,0), **@{-}; (52,0),**@{-};(65, 0), **@{-};
(0,0)*+{\bullet};(13,0)*+{\bullet};(26,0)*+{\bullet}; (39,0)*+{\bullet}; (52,0)*+{\bullet};
(0,-3)*+{p_1};(13,-3)*+{p_2};(26,-3)*+{p_3};(39,-3)*+{p_4};(52,-3)*+{p_5};(39,0)*+{\bullet};(52,0)*+{\bullet};
(0,0); (0,0); (-13,0), **@{-};
(26,0); (26,10), **@{-};(26,10)*+{\bullet};(23,10)*+{p_6};(52,0); (65,0), **@{-};(26,10); (26,20), **@{-};
(23,0)*+{\tiny{2t}};(29,0)*+{\tiny{2t}};(26,3)*+{\tiny{2t}};(35,0.6)*+{\tiny{-2t}};(42,0.6)*+{\tiny{4t}};
(48,0.6)*+{\tiny{-4t}};(55,0.6)*+{\tiny{6t}};(16,0.6)*+{\tiny{-2t}};
(10,0.6)*+{\tiny{4t}};(3,0.6)*+{\tiny{-4t}};(-2,0.6)*+{\tiny{6t}};(26,8)*+{\tiny{-2t}};(26,12)*+{\tiny{4t}};
\endxy
$$
Here a bullet corresponds to a fixed locus, an edge between two bullets represents the \emph{invariant} line linking the two fixed locus, and a ray represents an invariant direction. The corresponding weights of the torus action are given. Note that at the bullet $p_{3}$ there are three invariant directions, this means that $p_{3}$ is a \emph{fixed} $\mathbb{P}^{1}$. By the Atiyah-Bott localization formula,
\ben
\langle 1\rangle_{1,1,0}^{Y_{\hat{E}_{6}}}&=&2\int_{\overline{\cM}_{1,1}}\frac{(6t-\lambda_1)(-4t-\lambda_1)}{6t\cdot(-4t)}
+3\int_{\overline{\cM}_{1,1}}\frac{(4t-\lambda_1)(-2t-\lambda_1)}{4t\cdot(-2t)}\\
&&+\int_{\overline{\cM}_{1,1}\times \mathbb{P}^{1}}\frac{(2t-\lambda_1)(2H-\lambda_1)}{2t}\\
&=&\frac{1}{6t}\int_{\overline{\cM}_{1,1}}\lambda_1+\frac{3}{4t}\int_{\overline{\cM}_{1,1}}\lambda_1
-\frac{1}{t}\int_{\overline{\cM}_{1,1}}\lambda_1\\
&=&-\frac{1}{12t}\int_{\overline{\cM}_{1,1}}\lambda_1.
\een
where $2H=c_{1}(T_{\mathbb{P}^{1}})$.
\ben
\langle\cdot\rangle_{2,0,0}^{Y_{\hat{E}_{6}}}&=&2\int_{\overline{\cM}_{2,0}}\frac{((6t)^{2}-6t\lambda_1+\lambda_2)((-4t)^{2}-(-4t)\lambda_1
+\lambda_{2})}{6t\cdot(-4t)}
+3\int_{\overline{\cM}_{2,0}}\frac{((4t)^{2}-4t\lambda_1+\lambda_{2})((-2t)^{2}-(-2t)\lambda_1+\lambda_{2})}{4t\cdot(-2t)}\\
&&+\int_{\overline{\cM}_{2,0}\times \mathbb{P}^{1}}\frac{((2t)^{2}-2t\lambda_1+\lambda_2)(-2H\lambda_1+\lambda_{2})}{2t}\\
&=&\frac{1}{6t}\int_{\overline{\cM}_{2,0}}\lambda_1\lambda_2+\frac{3}{4t}\int_{\overline{\cM}_{2,0}}\lambda_1\lambda_2
-\frac{1}{t}\int_{\overline{\cM}_{2,0}}\lambda_1\lambda_{2}\\
&=&-\frac{1}{12t}\int_{\overline{\cM}_{2,0}}\lambda_1\lambda_2= -\frac{1}{2880t\cdot 24},
\een
where we use
\ben
\int_{\overline{\cM}_{2,0}}\lambda_1\lambda_2=\frac{1}{5760}
\een
from theorem 1 in \cite{FP}.\\
To verify that the other correlators vanishes, we follow the argument of \cite{BGh0707}. Denote the divisor  $[p_{1}p_{2}]$ defined by the line linking $p_{1}$ and $p_{2}$ by $E_{1}$, similarly $E_{2}=[p_{2}p_{3}]$, $E_{4}=[p_{3}p_{4}]$, $E_{5}=[p_{4}p_{5}]$, $E_{6}=[p_{3}p_{6}]$, and the divisor $[p_{3}]$ is denoted by $E_{3}$. The equivariant line bundles $L_{i}$ for $1\leq i\leq 6$ are defined such that $\alpha_{i}=c_{1}(L_{i})\}$  satisfies the intersection pairing
\bea\label{3}
\int_{E_{j}}\alpha_{i}=E_{i}\cdot E_{j},
\eea
which is the minus Cartan matrix by the classical McKay correnspondence. From the proof of \cite{BGh0707} we know that $L_{1}$ has weight $6t$ and $-2t$ at $p_{1}$ and $p_{2}$ respectively, and weight 0 at other  fixed loci. Thus
\ben
\langle \alpha_{1}\rangle_{1,1,0}^{Y_{\hat{E}_{6}}}&=&\int_{\overline{\cM}_{1,1}}6t\cdot\frac{(6t-\lambda_1)(-4t-\lambda_1)}{6t\cdot(-4t)}
+\int_{\overline{\cM}_{1,1}}(-2t)\cdot\frac{(4t-\lambda_1)(-2t-\lambda_1)}{4t\cdot(-2t)}=0.
\een
Similarly
$\langle \alpha_{5}\rangle_{1,1,0}^{Y_{\hat{E}_{6}}}=0$. $L_{2}$ has weight $4t$ at $p_2$, and weight zero at other fixed loci. Note also that $L_{2}\cong O(1)$ restricted to $p_{3}$, which easily follows from (\ref{3}). Thus
\ben
\langle \alpha_{2}\rangle_{1,1,0}^{Y_{\hat{E}_{6}}}&=&
\int_{\overline{\cM}_{1,1}}4t\cdot\frac{(4t-\lambda_1)(-2t-\lambda_1)}{4t\cdot(-2t)}+\int_{\overline{\cM}_{1,1}\times \mathbb{P}^{1}}H\cdot\frac{(2t-\lambda_1)(2H-\lambda_1)}{2t}=0.
\een
Similarly
$\langle \alpha_{4}\rangle_{1,1,0}^{Y_{\hat{E}_{6}}}=\langle \alpha_{6}\rangle_{1,1,0}^{Y_{\hat{E}_{6}}}=0$.  $L_{3}$ has weight $2t$ at $p_{3}$ and weight zero at other loci. Restricted to $p_{3}$, $L_{3}\cong O(-2)$. Thus
\ben
\langle \alpha_{3}\rangle_{1,1,0}^{Y_{\hat{E}_{6}}}&=&
\int_{\overline{\cM}_{1,1}\times \mathbb{P}^{1}}(2t-2H)\cdot\frac{(2t-\lambda_1)(2H-\lambda_1)}{2t}=0.
\een
Th  correlators for $(g,n,d)=(1,2,0)$, $(2,1,0)$ or $(3,0,0)$ are zero, because in these cases we need to integrate $\lambda_{i}^{2}$, $1\leq i\leq 3$, which is zero by Mumford's relations. The other correlators are zero obviously from the degree counting of the integrands and the string equation. Thus we complete the proof of the theorem for $G=\hat{E}_{6}$. In the following we give the proof of (\ref{5}) and (\ref{6}) for the other groups, and omit proof of  the vanishing statement which is similar to the case of $\hat{E}_{6}$.\\

$G=\hat{E}_{7}$:
$$ 
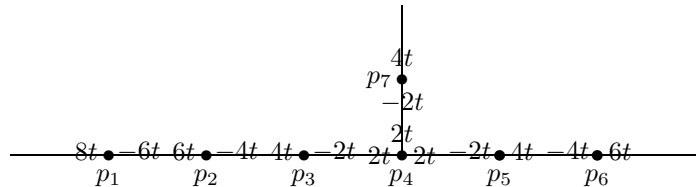
\xy
(0,0); (13,0), **@{-};(26,0), **@{-}; (26,0), **@{.}; (39,0),**@{-};(52, 0), **@{-};(-13, 0);(-26,0) **@{-};
(-13,0)*+{\bullet};(0,0)*+{\bullet};(13,0)*+{\bullet};(26,0)*+{\bullet}; (39,0)*+{\bullet}; (52,0)*+{\bullet};
(-13,-3)*+{p_1};(0,-3)*+{p_2};(13,-3)*+{p_3};(26,-3)*+{p_4};(39,-3)*+{p_5};(39,0)*+{\bullet};(52,0)*+{\bullet};(52,-3)*+{p_6};
(0,0); (0,0); (-13,0), **@{-};
(26,0); (26,10), **@{-};(26,10)*+{\bullet};(23,10)*+{p_7};(52,0); (65,0), **@{-};(26,10); (26,20), **@{-};
(23,0)*+{\tiny{2t}};(29,0)*+{\tiny{2t}};(26,3)*+{\tiny{2t}};(35,0.6)*+{\tiny{-2t}};(42,0,6)*+{\tiny{4t}};
(48,0.6)*+{\tiny{-4t}};(55,0.6)*+{\tiny{6t}};(17,0.6)*+{\tiny{-2t}};
(10,0.6)*+{\tiny{4t}};(4,0,6)*+{\tiny{-4t}};(-3,0,6)*+{\tiny{6t}};(26,7)*+{\tiny{-2t}};(26,13)*+{\tiny{4t}};(-9,0.6)*+{\tiny{-6t}};
(-16,0.6)*+{\tiny{8t}};
\endxy
$$

\ben
\langle 1\rangle_{1,1,0}^{Y_{\hat{E}_{7}}}&=&\int_{\overline{\cM}_{1,1}}\frac{(8t-\lambda_1)(-6t-\lambda_1)}{8t\cdot(-6t)}
+2\int_{\overline{\cM}_{1,1}}\frac{(6t-\lambda_1)(-4t-\lambda_1)}{6t\cdot(-4t)}\\
&&+3\int_{\overline{\cM}_{1,1}}\frac{(4t-\lambda_1)(-2t-\lambda_1)}{4t\cdot(-2t)}
+\int_{\overline{\cM}_{1,1}\times \mathbb{P}^{1}}\frac{(2t-\lambda_1)(2H-\lambda_1)}{2t}\\
&=&\frac{1}{24t}\int_{\overline{\cM}_{1,1}}\lambda_1+\frac{1}{6t}\int_{\overline{\cM}_{1,1}}\lambda_1+\frac{3}{4t}\int_{\overline{\cM}_{1,1}}\lambda_1
-\frac{1}{t}\int_{\overline{\cM}_{1,1}}\lambda_1\\
&=&-\frac{1}{24t}\int_{\overline{\cM}_{1,1}}\lambda_1=-\frac{1}{12t\cdot 48},
\een

\ben
\langle\cdot\rangle_{2,0,0}^{Y_{\hat{E}_{7}}}&=&\int_{\overline{\cM}_{2,0}}\frac{((8t)^{2}-8t\lambda_1+\lambda_2)((-6t)^{2}-(-6t)\lambda_1
+\lambda_{2})}{8t\cdot(-6t)}+2\int_{\overline{\cM}_{2,0}}\frac{((6t)^{2}-6t\lambda_1+\lambda_2)((-4t)^{2}-(-4t)\lambda_1
+\lambda_{2})}{6t\cdot(-4t)}\\
&&+3\int_{\overline{\cM}_{2,0}}\frac{((4t)^{2}-4t\lambda_1+\lambda_{2})((-2t)^{2}-(-2t)\lambda_1+\lambda_{2})}{4t\cdot(-2t)}
+\int_{\overline{\cM}_{2,0}\times \mathbb{P}^{1}}\frac{((2t)^{2}-2t\lambda_1+\lambda_2)(-2H\lambda_1+\lambda_{2})}{2t}\\
&=&\frac{1}{24t}\int_{\overline{\cM}_{2,0}}\lambda_1\lambda_2+\frac{1}{6t}\int_{\overline{\cM}_{2,0}}\lambda_1\lambda_2+\frac{3}{4t}\int_{\overline{\cM}_{2,0}}\lambda_1\lambda_2
-\frac{1}{t}\int_{\overline{\cM}_{2,0}}\lambda_1\lambda_{2}\\
&=&-\frac{1}{24t}\int_{\overline{\cM}_{2,0}}\lambda_1\lambda_2= -\frac{1}{2880t\cdot 48}.
\een

$G=\hat{E}_{8}$:

$$ 
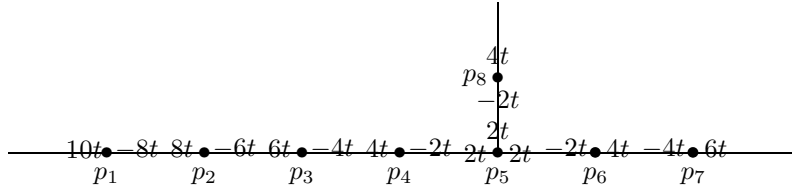
\xy
(0,0); (13,0), **@{-};(26,0), **@{-}; (26,0), **@{.}; (39,0),**@{-};(52, 0), **@{-};(-13, 0);(-26,0) **@{-};
(-13,0)*+{\bullet};(0,0)*+{\bullet};(13,0)*+{\bullet};(26,0)*+{\bullet}; (39,0)*+{\bullet}; (52,0)*+{\bullet};
(-13,-3)*+{p_2};(0,-3)*+{p_3};(13,-3)*+{p_4};(26,-3)*+{p_5};(39,-3)*+{p_6};(39,0)*+{\bullet};(52,0)*+{\bullet};(52,-3)*+{p_7};
(0,0); (0,0); (-13,0), **@{-};
(26,0); (26,10), **@{-};(26,10)*+{\bullet};(23,10)*+{p_8};(52,0); (65,0), **@{-};(26,10); (26,20), **@{-};
(23,0)*+{\tiny{2t}};(29,0)*+{\tiny{2t}};(26,3)*+{\tiny{2t}};(35,0.6)*+{\tiny{-2t}};(42,0,6)*+{\tiny{4t}};
(48,0.6)*+{\tiny{-4t}};(55,0.6)*+{\tiny{6t}};(17,0.6)*+{\tiny{-2t}};
(10,0.6)*+{\tiny{4t}};(4,0,6)*+{\tiny{-4t}};(-3,0,6)*+{\tiny{6t}};(26,7)*+{\tiny{-2t}};(26,13)*+{\tiny{4t}};(-9,0.6)*+{\tiny{-6t}};
(-16,0.6)*+{\tiny{8t}};
(-26,0)*+{\bullet}; (-26,-3)*+{p_{1}}; (-29,0.5)*{\tiny{10t}};(-22,0.5)*+{\tiny{-8t}};
(-26,0); (-39,0),**@{-};
\endxy
$$

\ben
\langle 1\rangle_{1,1,0}^{Y_{\hat{E}_{8}}}&=&\int_{\overline{\cM}_{1,1}}\frac{(10t-\lambda_1)(-8t-\lambda_1)}{10t\cdot(-8t)}
+\int_{\overline{\cM}_{1,1}}\frac{(8t-\lambda_1)(-6t-\lambda_1)}{8t\cdot(-6t)}
+2\int_{\overline{\cM}_{1,1}}\frac{(6t-\lambda_1)(-4t-\lambda_1)}{6t\cdot(-4t)}\\
&&+3\int_{\overline{\cM}_{1,1}}\frac{(4t-\lambda_1)(-2t-\lambda_1)}{4t\cdot(-2t)}
+\int_{\overline{\cM}_{1,1}\times \mathbb{P}^{1}}\frac{(2t-\lambda_1)(2H-\lambda_1)}{2t}\\
&=&\frac{1}{40t}\int_{\overline{\cM}_{1,1}}\lambda_1+\frac{1}{24t}\int_{\overline{\cM}_{1,1}}\lambda_1+\frac{1}{6t}\int_{\overline{\cM}_{1,1}}\lambda_1+\frac{3}{4t}\int_{\overline{\cM}_{1,1}}\lambda_1
-\frac{1}{t}\int_{\overline{\cM}_{1,1}}\lambda_1\\
&=&-\frac{1}{60t}\int_{\overline{\cM}_{1,1}}\lambda_1=-\frac{1}{12t\cdot 120},
\een

\ben
\langle\cdot\rangle_{2,0,0}^{Y_{\hat{E}_{8}}}&=&\int_{\overline{\cM}_{2,0}}\frac{((10t)^{2}-10t\lambda_1+\lambda_2)((-8t)^{2}-(-8t)\lambda_1
+\lambda_{2})}{10t\cdot(-8t)}+\int_{\overline{\cM}_{2,0}}\frac{((8t)^{2}-8t\lambda_1+\lambda_2)((-6t)^{2}-(-6t)\lambda_1
+\lambda_{2})}{8t\cdot(-6t)}\\
&&+2\int_{\overline{\cM}_{2,0}}\frac{((6t)^{2}-6t\lambda_1+\lambda_2)((-4t)^{2}-(-4t)\lambda_1
+\lambda_{2})}{6t\cdot(-4t)}
+3\int_{\overline{\cM}_{2,0}}\frac{((4t)^{2}-4t\lambda_1+\lambda_{2})((-2t)^{2}-(-2t)\lambda_1+\lambda_{2})}{4t\cdot(-2t)}\\
&&+\int_{\overline{\cM}_{2,0}\times \mathbb{P}^{1}}\frac{((2t)^{2}-2t\lambda_1+\lambda_2)(-2H\lambda_1+\lambda_{2})}{2t}\\
&=&\frac{1}{40t}\int_{\overline{\cM}_{2,0}}\lambda_1\lambda_2+\frac{1}{24t}\int_{\overline{\cM}_{2,0}}\lambda_1\lambda_2+\frac{1}{6t}\int_{\overline{\cM}_{2,0}}\lambda_1\lambda_2+\frac{3}{4t}\int_{\overline{\cM}_{2,0}}\lambda_1\lambda_2
-\frac{1}{t}\int_{\overline{\cM}_{2,0}}\lambda_1\lambda_{2}\\
&=&-\frac{1}{60t}\int_{\overline{\cM}_{2,0}}\lambda_1\lambda_2= -\frac{1}{2880t\cdot 120}.
\een

$G=\hat{A}_{n}$:
$$ \xy
(53,1.7)*+{\tiny{(n+1)t_{1}}};(68,1.7)*+{\tiny{-nt_{1}+t_{2}}};(60,-3)*+{p_{1}};
(45,0);(60,0),**@{-}; (60,0);(75,0),**@{-};(60,0)*+{\bullet};(75,0);(105,0),**@{.};(105,0);(120,0),**@{-};
(100,1.7)*+{\tiny{(n+2-k)t_{1}-(k-1)t_{2}}};(138,1.7)*+{\tiny{-(n+1-k)t_{1}+kt_{2}}};(120,-3)*+{p_{k}};
(120,0)*+{\bullet};(120,0);(135,0),**@{-};(135,0);(165,0),**@{.};(165,0);(180,0),**@{-};(180,0)*+{\bullet};
(180,0);(195,0),**@{-}; (180,-3)*+{p_{n+1}};(173,1.7)*+{\tiny{t_{1}-nt_{2}}};(188,1.7)*+{\tiny{(n+1)t_{2}}};
\endxy
$$

\ben
\langle 1\rangle_{1,1,0}^{Y_{\hat{A}_{n}}}&=&\sum_{k=1}^{n+1}\int_{\overline{\cM}_{1,1}}\frac{((n+2-k)t_{1}-(k-1)t_{2}-\lambda_1)
(-(n+1-k)t_{1}+k t_{2}-\lambda_1)}{\big(((n+2-k)t_{1}-(k-1)t_{2}\big)\big(-(n+1-k)t_{1}+k t_{2}\big)}\\
&=& \sum_{k=1}^{n+1}\Big(-\frac{1}{(n+2-k)t_{1}-(k-1)t_{2}}+\frac{1}{(n+1-k)t_{1}-k t_{2}}\Big)\int_{\overline{\cM}_{1,1}}\lambda_{1}\\
&=& -\frac{t_{1}+t_{2}}{24(n+1)t_{1}t_{2}},
\een

\ben
\langle\cdot\rangle_{2,0,0}^{Y_{\hat{A}_{n}}}&=&\sum_{k=1}^{n+1}\int_{\overline{\cM}_{2,0}}\frac{\big((n+2-k)t_{1}-(k-1)t_{2}\big)^{2}
-\big((n+2-k)t_{1}-(k-1)t_{2}\big)\lambda_1+\lambda_2}{(n+2-k)t_{1}-(k-1)t_{2}}\\
&&\cdot\frac{\big(-(n+1-k)t_{1}+k t_{2}\big)^{2}-\big(-(n+1-k)t_{1}+k t_{2}\big)\lambda_1
+\lambda_{2}}{-(n+1-k)t_{1}+k t_{2}}\\
&=& \sum_{k=1}^{n+1}\Big(-\frac{1}{(n+2-k)t_{1}-(k-1)t_{2}}+\frac{1}{(n+1-k)t_{1}-k t_{2}}\Big)\int_{\overline{\cM}_{2,0}}\lambda_{1}\lambda_{2}\\
&=& -\frac{t_{1}+t_{2}}{5760(n+1)t_{1}t_{2}}.
\een

$G=\hat{D}_{n}$:
$$ \xy
(52,1.7)*+{\tiny{(2n-4) t}};(70,1.7)*+{\tiny{(-2n+6)t}};(60,-3)*+{p_{1}};
(45,0);(60,0),**@{-}; (60,0);(75,0),**@{-};(60,0)*+{\bullet};(75,0);(105,0),**@{.};(105,0);(120,0),**@{-};
(108,1.7)*+{\tiny{(2n-2k-2)t}};(133,1.7)*+{\tiny{(-2n+2k+4)t}};(120,-3)*+{p_{k}};
(120,0)*+{\bullet};(120,0);(135,0),**@{-};(135,0);(165,0),**@{.};(165,0);(180,0),**@{-};(180,0)*+{\bullet};
(180,0);(195,10),**@{-}; (180,-3)*+{p_{n-2}};(177,0)*+{2t};
(180,0);(195,-10),**@{-};(195,10)*+{\bullet};(195,-10)*+{\bullet};(183,3)*+{2t};(183,-1)*+{2t};
(192,7)*+{-2t};(192,-7)*+{-2t};
(195,10);(205,10),**@{-};(195,-10);(205,-10),**@{-};(198,10.6)*+{4t};(198,-9.4)*+{4t};
\endxy
$$

\ben
\langle 1\rangle_{1,1,0}^{Y_{\hat{D}_{n}}}&=&\sum_{k=1}^{n-3}\int_{\overline{\cM}_{1,1}}\frac{\big((2n-2k-2)t-\lambda_1\big)
\big((-2n+2k+4)t-\lambda_1\big)}{(2n-2k-2)(-2n+2k+4)t^{2}}\\
&&+2\int_{\overline{\cM}_{1,1}}\frac{(4t-\lambda_1)(-2t-\lambda_1)}{4t\cdot(-2t)}+\int_{\overline{\cM}_{1,1}\times \mathbb{P}^{1}}\frac{(2t-\lambda_1)(2H-\lambda_1)}{2t}\\
&=& \sum_{k=1}^{n-3}\Big(-\frac{1}{(2n-2k-2)t}+\frac{1}{(2n-2k-4)t}\Big)\int_{\overline{\cM}_{1,1}}\lambda_{1}+
\frac{1}{2t}\int_{\overline{\cM}_{1,1}}\lambda_1
-\frac{1}{t}\int_{\overline{\cM}_{1,1}}\lambda_1\\
&=& -\frac{1}{(2n-4)t\cdot 24}=-\frac{1}{(4n-8)\cdot 12t},
\een

\ben
\langle\cdot\rangle_{2,0,0}^{Y_{\hat{D}_{n}}}&=&\sum_{k=1}^{n-3}\int_{\overline{\cM}_{2,0}}\frac{\big((2n-2k-2)^{2}t^{2}
-(2n-2k-2)t\lambda_1+\lambda_2\big)\big((-2n+2k+4)^{2}t^{2}
-(-2n+2k+4)t\lambda_1+\lambda_2\big)}{(2n-2k-2)(-2n+2k+4)t^{2}}\\
&&+2\int_{\overline{\cM}_{2,0}}\frac{((4t)^{2}-4t\lambda_1+\lambda_{2})((-2t)^{2}-(-2t)\lambda_1+\lambda_{2})}{4t\cdot(-2t)}
+\int_{\overline{\cM}_{2,0}\times \mathbb{P}^{1}}\frac{((2t)^{2}-2t\lambda_1+\lambda_2)(-2H\lambda_1+\lambda_{2})}{2t}\\
&=& -\frac{1}{(4n-8)\cdot 2880t}.
\een
\hfill\qedsymbol\\

\section{Equivariant Gromov-Witten invariants of $[\mathbb{C}^2/G]$}
For a finite group G, consider the following diagram with a cartesian square

\ben
\xymatrix{ \mathcal{V}\ar[r]^{g}\ar[d]^{p} & pt \ar[d] \\
\mathcal{U} \ar[r]^{f}\ar[d]^{\pi}& BG\\
\overline{\cM}_{g,n}(BG; \gamma_{1},\cdots,\gamma_{n})& }
\een

When $(\gamma_{1},\cdots,\gamma_{n})=(\lbr 1\rbr,\cdots,\lbr 1\rbr)$, there is a component in $\overline{\cM}_{g,n}(BG; \gamma_{1},\cdots,\gamma_{n})$ for which $p$ is a trivial $G$-torsor. We denote this component by $\overline{\cM}_{g,n}^{\text{trivial}}(BG;\lbr 1\rbr,\cdots,\lbr 1\rbr)$, and the remaining component by $\overline{\cM}_{g,n}^{\text{nontri}}(BG;\lbr 1\rbr,\cdots,\lbr 1\rbr)$. \\

\emph{Assume from now on that $G$ is a finite subgroup of $SL(2,\mathbb{C})$}, and denote the standard 2-dimensional representation induced by this inclusion by $V_{\rho_{1}}$. The corresponding vector bundle on $BG$ is also denote by $V_{\rho_{1}}$, by an abuse of notation. Thus we have
\begin{proposition}\label{15}
Suppose the number of $\lbr 1\rbr$ in $\gamma_{1},\cdots,\gamma_{n}$ is $m$. Then\\
(i) $R^{0}\pi_{*}f^{*}V_{\rho_{1}}$ and $R^{1}\pi_{*}f^{*}V_{\rho_{1}}$ are vector bundles on every components;\\
(ii)For $G$ of type $A$ (resp., type $D$ or type $E$), $R^{0}\pi_{*}f^{*}V_{\rho_{1}}$ is a trivial rank $2$ (resp. $1$) vector bundle on $\overline{\cM}_{g,n}^{\text{trivial}}(BG;\lbr 1\rbr,\cdots,\lbr 1\rbr)$, and is rank $0$ on all the other components;\\
(iii) For $G$ of type $A$ (resp., type $D$ or type $E$), the rank of $R^{1}\pi_{*}f^{*}V_{\rho_{1}}$  on $\overline{\cM}_{g,n}^{\text{trivial}}(BG;\lbr 1\rbr,\cdots,\lbr 1\rbr)$ is $2g$ (resp. $2g-1$), on the other components is $2g-2+n-m$.
\end{proposition}
\emph{Proof}: The statement (i) follows from (ii). The statement (iii) follows from (ii) by the orbifold Riemann-Roch theorem. Thus we need only to check (ii). Let $\mathcal{C}$ be a fibre of $\pi$, and $\tilde{\mathcal{C}}$ be the corresponding $G$-torsor over $\mathcal{C}$. We have $H^{0}(f^{*}V_{\rho_{1}})=(H^{0}(\mathcal{O}_{\tilde{\mathcal{C}}})\otimes V_{\rho_{1}})^{G}$. From the Dynkin diagram of the ADE-singularities, we easily see that, for an irreducible representation $V$ of $G$, the existence of at least one trivial summand in the decomposition of $V\otimes V_{\rho_{1}}$ implies that $V$ is a faithful representation. Thus $(H^{0}(\mathcal{O}_{\tilde{\mathcal{C}}})\otimes V_{\rho_{1}})^{G}\neq 0$ forces the $G$-torsor $\tilde{\mathcal{C}}\rightarrow\mathcal{C} $ to be trivial; and when $\tilde{\mathcal{C}}\rightarrow\mathcal{C} $ is trivial, the number of trivial summands in $H^{0}(\mathcal{O}_{\tilde{\mathcal{C}}})\otimes V_{\rho_{1}}$ reads from that of $\mathbb{C}G\otimes V_{\rho_{1}}$, where $\mathbb{C}G$ denotes the regular representation of the finite group $G$. It is well known that the irreducible summands of $\mathbb{C}G$ runs over each irreducible representation of $G$ for one time, thus again from the Dynkin diagram we see that, for $G$ of type $A$ (resp., type $D$ or type $E$) the number of  trivial summands in $H^{0}(\mathcal{O}_{\tilde{\mathcal{C}}})\otimes V_{\rho_{1}}$ is two (resp., one). The triviality of the corresponding vector bundle is straightforward. \hfill\qedsymbol\\

Therefore to compute nonzero primary orbifold Gromov-Witten invariants of $[\mathbb{C}^{2}/G]$, there are only three cases we need to consider:\\
(i)$g=0$;\\
(ii)For $G$ of type $A$ (resp., type $D$ or type $E$), $g+n\leq 3$ (resp., $g+n\leq 2$ ), $\overline{\cM}_{g,n}^{\text{trivial}}(BG;\lbr 1\rbr,\cdots,\lbr 1\rbr)$;\\
(iii)$g=1$, $m=0$, $\overline{\cM}_{g,n}^{\text{nontri}}(BG;\lbr 1\rbr,\cdots,\lbr 1\rbr)$. \\

The case (i) has been stated as
the quantum McKay correspondence conjecture \ref{13}. For case (ii), we have
\ben
&&\int_{\overline{\cM}_{1,1}^{\text{trivial}}(BG;\lbr 1\rbr)}\frac{(t-\lambda_{1})(t-\lambda_{1})}{t^{2}}\\
&=&-\frac{2}{t}\frac{1}{|G|}\int_{\overline{\cM}_{1,1}}\lambda_{1}\\
&=&-\frac{1}{12|G|t},
\een

\ben
&&\int_{\overline{\cM}_{2,0}^{\text{trivial}}(BG)}\frac{(t^{2}-t\lambda_{1}+\lambda_{2})(t^{2}-t\lambda_{1}+\lambda_{2})}{t^2}\\
&=&-\frac{2}{t}\frac{1}{|G|}\int_{\overline{\cM}_{2,0}}\lambda_{1}\lambda_{2}\\
&=&-\frac{1}{2880t|G|},
\een
and other integrals vanish by the Mumford relation. By theorem \ref{4}, these computations establish the correspondence between the degree zero invariants on the resolution spaces on the Hurwitz-Hodge integrals on $\overline{\cM}_{g,n}^{\text{trivial}}(BG;\lbr 1\rbr,\cdots,\lbr 1\rbr)$ in higher genera , for ADE singularities. For ancestor invariants similar computations show that the same correspondence holds.\\

In case (iii) we make the
\begin{Conjecture}\label{7}
\bea\label{1}
\langle e_{ \gamma_{1}}, \cdots, e_{ \gamma_{n}}\rangle^{[\mathbb{C}^{2}/G]}_{1,n}=0
\eea
 for $(  \gamma_{1},\cdots,  \gamma_{n})\neq (\lbr 1\rbr,\cdots,\lbr 1\rbr)$ and $n\geq 1$.
\end{Conjecture}

We check this conjecture for $n=1$ by the orbifold quantum Riemann-Roch theorem (\cite{Ts}),
\begin{proposition}\label{32}
\bea\label{2}
\langle e_{ \gamma}\rangle^{[\mathbb{C}^{2}/G]}_{1,1}=0
\eea
 for $\gamma\neq \lbr 1\rbr$.
\end{proposition}
The proof is case by case checking, which is somewhat tedious and we give it in the appendix A.

By a comparison to theorem \ref{4}, we obtain
\begin{theorem}\label{8}
Assuming the quantum McKay correspondence conjecture \ref{13}, and the conjecture \ref{7}, the Gromov-Witten potential $F_{g}^{\widehat{\mathbb{C}^{2}/G}}$ (restricted to the absolute stable range), after the analytic continuation and the change of variables \ref{12}, equals the Gromov-Witten potential $F_{g}^{[\mathbb{C}^{2}/G]}$, for every $g\geq 0$.
\end{theorem}\hfill\qedsymbol

\section{The crepant resolution conjecture for gravitational ancestors}
In this section we give an ancestor version of theorem \ref{8}. First as an analogy to conjecture \ref{7}, we make the following
\begin{Conjecture}\label{9}
Let $\rho$ be an irreducible 2-dim representation of a finite group $G$ belonging to ADE types. Then
\ben
\int_{\overline{\cM}_{g,n}(BG;\lbr a_1\rbr,\cdots,\lbr a_n\rbr)}c_{2g-2+n}(F_{\rho}^{1})\prod_{k=1}^{n}\bar{\psi}_{\lbr a_k\rbr}^{l_k}=0,
\een
when $(\lbr a_1\rbr,\cdots,\lbr a_n\rbr)\neq \lbr 1\rbr^n$, and $\sum_{k=1}^{n}l_{k}=g-1$.
\end{Conjecture}

Conjecture \ref{7} is a special case of this one. Although $[\mathbb{C}^{2}/G]$ has a canonical \emph{orbifold holomorphic symplectic form}, the argument of \cite{OP} seems not able to directly extend to this conjecture.\\

When conjecture \ref{9} holds, we can apply the method of \cite{FSZ} to prove an ancestor version of theorem \ref{8}.

\begin{theorem}\label{10}
Assuming the quantum McKay correspondence \ref{13}, and the conjecture \ref{9}, the total ancestor Gromov-Witten potential $\mathcal{A}^{\widehat{\mathbb{C}^{2}/G}}$, after the analytic continuation from and the change of variables \ref{12}, equals the total ancestor Gromov-Witten potential $\mathcal{A}^{[\mathbb{C}^{2}/G]}$.
\end{theorem}
\emph{Proof}: By proposition \ref{15}, and assuming  the conjecture \ref{9}, every nonzero ancestor invariant either has $\bar{\psi}$ classes of total degree at least $g$, or has been treated in the discussions following proposition \ref{15}. Thus as in the proof of theorem 5 in \cite{FSZ}, the ancestor correlators can be reduced to the primary correlators by the tautological relations. Note also that every tautological relations can be expressed as a partial differential equation, which is \emph{independent of $q_{i}$'s}. Thus the theorem follows from theorem \ref{8}.\hfill\qedsymbol\\

By the proposition 3.2 of \cite{BGP}, conjecture \ref{9} holds for groups of type A, thus we have
\begin{corollary}\label{14}
 The total ancestor Gromov-Witten potential $\mathcal{A}^{\widehat{\mathbb{C}^{2}/\mathbb{Z}_{n}}}$, after the analytic continuation  and the change of variables \ref{12}, equals the total ancestor Gromov-Witten potential $\mathcal{A}^{[\mathbb{C}^{2}/\mathbb{Z}_{n}]}$.
\end{corollary}\hfill\qedsymbol\\

\begin{remark}\label{14}
By the comparison formula between ancestor and descendant invariants \cite{KM}, it's not hard to see that the coefficents of the equivariant parameter $t$ in the stationary part of the total ancestor potential is equal to that of the total descendant potential. Thus in this way we recover the higher genera part of the corresponding result in \cite{Zhou0811} (but we \emph{need} to use the genus zero part of it!).
\end{remark}

    \begin{appendices}

\section{Proof of the proposition \ref{32}}
In this appendix we verify the proposition \ref{32}. The Hurwitz-Hodge bundle $R^{1}\pi_{*}f^{*}V_{\rho_{1}}$ on $\overline{\cM}_{1,1}(BG;\lbr\gamma\rbr)$ is of rank 1 for $G$ of ADE types and $\rho_{1}$ the standard 2-dim representation and a nontrivial conjugacy class $\gamma$. Thus what we need to prove amounts to
\bea
\langle e_{\lbr \gamma\rbr}ch_{1}^{\rho_{1}}\rangle_{1,1}^{BG}=0,
\eea
where we use $ch_{1}^{\rho_{1}}$ to denote the first component of the chern character of $R^{1}\pi_{*}f^{*}V_{\rho_{1}}$ for short.
Our strategy is to directly use the quantum Riemann-Roch formula \cite{Ts} to compute this integral (see also \cite{Zhou0710}).\\

 We shall give  the details for groups of type $D$ and $E_{8}$. For $E_{7}$ we only give the length three correlators and the differential operators. $E_{6}$ is a normal subgroup of $E_{7}$, and  the corresponding length three correlators and the differential operators are easily deduced from those of $E_{7}$, thus we omit them. All the information about the generators and  the relations of the groups, as well the table of characters, are borrowed from the appendix A of \cite{Ste}.\\
\subsection{The binary dihedral groups}
Let $\mathbf{c}(\cdot)=\exp(\sum_{k=1}s_{k}ch_{k}^{\rho_{1}})$, where $(s_{1},\cdots)$ is series of formal parameters. We refer the readers to \cite{Ts} for the definition of $\mathbf{c}(\cdot)$-twisted Gromov-Witten invariants of $\mathcal{B}\hat{D}_{n}$ and Givental's quadratic quantization formalism (see also \cite{Zhou0710}). For example
\ben
F_1^{tw}=\sum_{n=1}^{\infty}\frac{1}{n!}\int_{\overline{\cM}_{1,n}(B\hat{D}_{m})}\exp(\sum_{s\geqslant 1}s_k\ch_{k} )\prod_{i=1}^{n}\sum_{l=1}^{\infty}\ev_{i}^{\ast}(\sum_{\llbracket\gamma\rrbracket}t_{l}^{\llbracket\gamma\rrbracket}
e_{\llbracket\gamma\rrbracket})\bar{\psi}_{\llbracket\gamma\rrbracket}^{l},
\een
where we have abbreviate the superscript $^{\rho_{1}}$ in $ch_{k}^{\rho_{1}})$ since no other characters will be considered in this subsection.
The quadratic operators  are
\ben
&&\big(\frac{A_{p+1}(V_{\rho_{1}})z^p}{(p+1)!}\big)^{\wedge}\\
&=&\frac{2B_{p+1}}{(p+1)!}\pd_{\lbr 1\rbr,1+p}
-\frac{2B_{p+1}}{(p+1)!}\sum_{l= 0}^{\infty}t_{l}^{\lbr 1\rbr}\pd_{\lbr 1\rbr,l+p}-\sum_{k=1}^{n-3}\frac{B_{p+1}(\frac{k}{2n-4})+B_{p+1}(\frac{2n-4-k}{2n-4})}{(p+1)!}\sum_{l= 0}^{\infty}t_{l}^{\lbr a^k\rbr}\pd_{\lbr a^k\rbr,l+p}\\
&&-\frac{2B_{p+1}(\frac{n-2}{2n-4})}{(p+1)!}\sum_{l=0}^{\infty}t_{l}^{\lbr a^{n-2}\rbr}\pd_{\lbr a^{n-2}\rbr,l+p}-\frac{B_{p+1}(\frac{1}{4})+B_{p+1}(\frac{3}{4})}{(p+1)!}\sum_{l=0}^{\infty}(t_{l}^{\lbr b\rbr}\pd_{\lbr b\rbr,l+p}+t_{l}^{\lbr ab\rbr}\pd_{\lbr ab\rbr,l+p})
\\
&&+\frac{\hbar^2}{2}\sum_{l= 0}^{p-1}(-1)^l\Bigg((4n-8)\frac{2B_{p+1}}{(p+1)!}\pd_{\lbr 1\rbr,l}\pd_{\lbr 1\rbr,p-1-l}+(2n-4)\sum_{k=1}^{n-3}\frac{B_{p+1}(\frac{k}{2n-4})+B_{p+1}(\frac{2n-4-k}{2n-4})}{(p+1)!}\\
&& \cdot\pd_{\lbr a^k \rbr,l}\pd_{\lbr a^k\rbr,p-1-l}+(4n-8)\frac{2B_{p+1}(\frac{1}{2})}{(p+1)!}\pd_{\lbr a^{n-2}\rbr,l}\pd_{\lbr a^{n-2}\rbr,p-1-l}+
4\cdot\frac{B_{p+1}(\frac{1}{4})+B_{p+1}(\frac{3}{4})}{(p+1)!}\\
&&\cdot(\pd_{\lbr b\rbr,l}\pd_{
\lbr b \rbr,p-1-l}+\pd_{\lbr ab\rbr,l}\pd_{
\lbr ab \rbr,p-1-l})\Bigg).
\een
Then by the orbifold quantum Riemann-Roch theorem (\cite{Ts}), we have
\ben
&&\sum_{n=1}^{\infty}\frac{1}{n!}\int_{\overline{\cM}_{1,n}(B\hat{D}_{m})}\ch_{p}\exp(\sum_{s\geqslant 1}s_k\ch_{k} )\prod_{i=1}^{n}\sum_{k=1}^{\infty}\ev_{i}^{\ast}(\sum_{\llbracket\gamma\rrbracket}t_{k}^{\llbracket\gamma\rrbracket}
e_{\llbracket\gamma\rrbracket})\bar{\psi}_{\llbracket\gamma\rrbracket}\\
&=&\frac{2B_{p+1}}{(p+1)!}\pd_{\lbr 1\rbr,1+p}F_{1}^{tw}
-\frac{2B_{p+1}}{(p+1)!}\sum_{l= 0}^{\infty}t_{l}^{\lbr 1\rbr}\pd_{\lbr 1\rbr,l+p}F_{1}^{tw}-\sum_{k=1}^{n-3}\frac{B_{p+1}(\frac{k}{2n-4})+B_{p+1}(\frac{2n-4-k}{2n-4})}{(p+1)!}\sum_{l= 0}^{\infty}t_{l}^{\lbr a^k\rbr}\pd_{\lbr a^k\rbr,l+p}F_{1}^{tw}\\
&&-\frac{2B_{p+1}(\frac{n-2}{2n-4})}{(p+1)!}\sum_{l=0}^{\infty}t_{l}^{\lbr a^{n-2}\rbr}\pd_{\lbr a^{n-2}\rbr,l+p}F_{1}^{tw}-\frac{B_{p+1}(\frac{1}{4})+B_{p+1}(\frac{3}{4})}{(p+1)!}\sum_{l=0}^{\infty}(t_{l}^{\lbr b\rbr}\pd_{\lbr b\rbr,l+p}+t_{l}^{\lbr ab\rbr}\pd_{\lbr ab\rbr,l+p})F_{1}^{tw}
\\
&&+\frac{1}{2}\sum_{l= 0}^{p-1}(-1)^l\Bigg((4n-8)\frac{2B_{p+1}}{(p+1)!}\pd_{\lbr 1\rbr,l}\pd_{\lbr 1\rbr,p-1-l}+(2n-4)\sum_{k=1}^{n-3}\frac{B_{p+1}(\frac{k}{2n-4})+B_{p+1}(\frac{2n-4-k}{2n-4})}{(p+1)!}\\
&& \cdot\pd_{\lbr a^k \rbr,l}\pd_{\lbr a^k\rbr,p-1-l}+(4n-8)\frac{2B_{p+1}(\frac{1}{2})}{(p+1)!}\pd_{\lbr a^{n-2}\rbr,l}\pd_{\lbr a^{n-2}\rbr,p-1-l}+
4\cdot\frac{B_{p+1}(\frac{1}{4})+B_{p+1}(\frac{3}{4})}{(p+1)!}\\
&&\cdot(\pd_{\lbr b\rbr,l}\pd_{
\lbr b \rbr,p-1-l}+\pd_{\lbr ab\rbr,l}\pd_{
\lbr ab \rbr,p-1-l})\Bigg)F_{0}^{tw}\\
&&+\frac{1}{2}\sum_{l= 0}^{p-1}(-1)^l\Bigg((4n-8)\frac{2B_{p+1}}{(p+1)!}\pd_{\lbr 1\rbr,l}F_{0}^{tw}\pd_{\lbr 1\rbr,p-1-l}F_{1}^{tw}+(2n-4)\sum_{k=1}^{n-3}\frac{B_{p+1}(\frac{k}{2n-4})+B_{p+1}(\frac{2n-4-k}{2n-4})}{(p+1)!}\\
&& \cdot\pd_{\lbr a^k \rbr,l}F_{0}^{tw}\pd_{\lbr a^k\rbr,p-1-l}F_{1}^{tw}+(4n-8)\frac{2B_{p+1}(\frac{1}{2})}{(p+1)!}\pd_{\lbr a^{n-2}\rbr,l}F_{0}^{tw}\pd_{\lbr a^{n-2}\rbr,p-1-l}F_{1}^{tw}+
4\cdot\frac{B_{p+1}(\frac{1}{4})+B_{p+1}(\frac{3}{4})}{(p+1)!}\\
&&\cdot(\pd_{\lbr b\rbr,l}F_{0}^{tw}\pd_{
\lbr b \rbr,p-1-l}F_{1}^{tw}+\pd_{\lbr ab\rbr,l}F_{0}^{tw}\pd_{
\lbr ab \rbr,p-1-l}F_{1}^{tw})\Bigg)\\
&&+\frac{1}{2}\sum_{l= 0}^{p-1}(-1)^l\Bigg((4n-8)\frac{2B_{p+1}}{(p+1)!}\pd_{\lbr 1\rbr,l}F_{1}^{tw}\pd_{\lbr 1\rbr,p-1-l}F_{0}^{tw}+(2n-4)\sum_{k=1}^{n-3}\frac{B_{p+1}(\frac{k}{2n-4})+B_{p+1}(\frac{2n-4-k}{2n-4})}{(p+1)!}\\
&& \cdot\pd_{\lbr a^k \rbr,l}F_{1}^{tw}\pd_{\lbr a^k\rbr,p-1-l}F_{0}^{tw}+(4n-8)\frac{2B_{p+1}(\frac{1}{2})}{(p+1)!}\pd_{\lbr a^{n-2}\rbr,l}F_{1}^{tw}\pd_{\lbr a^{n-2}\rbr,p-1-l}F_{0}^{tw}+
4\cdot\frac{B_{p+1}(\frac{1}{4})+B_{p+1}(\frac{3}{4})}{(p+1)!}\\
&&\cdot(\pd_{\lbr b\rbr,l}F_{1}^{tw}\pd_{
\lbr b \rbr,p-1-l}F_{0}^{tw}+\pd_{\lbr ab\rbr,l}F_{1}^{tw}\pd_{
\lbr ab \rbr,p-1-l}F_{0}^{tw})\Bigg).\een

In the following computations, we need the values of the genus zero length three correlators, which have been given in \cite{Hu}. Thus for $1\leq k < \frac{n-2}{2}$,
\bea\label{31}
&&\langle ch_{1}e_{\lbr a^{2k}\rbr}\rangle_{1,1}^{\mathcal{B}\hat{D}_{n}}\nonumber\\
&=&B_{2}\langle e_{\lbr a^{2k}\rbr}(e_{\lbr 1\rbr}\bar{\psi}_{2}^{2})\rangle_{1,2}^{\mathcal{B}\hat{D}_{n}}-\frac{B_{2}(\frac{2k}{2n-4})+B_{2}(\frac{2n-4-2k}{2n-4})}{2!}
\langle e_{\lbr a^{2k}\rbr}\bar{\psi}_{1}^{1}\rangle_{1,1}^{\mathcal{B}\hat{D}_{n}}\nonumber\\
&&+(n-2)\frac{B_{2}(\frac{n-2-k}{2n-4})+B_{2}(\frac{n-2+k}{2n-4})}{2!}\langle e_{\lbr a^{2k}\rbr}e_{\lbr a^{n-2-k}\rbr}e_{\lbr a^{n-2-k}\rbr}\rangle_{0,3}^{\mathcal{B}\hat{D}_{n}}\nonumber\\
&&+(n-2)\frac{B_{2}(\frac{k}{2n-4})+B_{2}(\frac{2n-4-k}{2n-4})}{2!}\langle e_{\lbr a^{2k}\rbr}e_{\lbr a^{k}\rbr}e_{\lbr a^{k}\rbr}\rangle_{0,3}^{\mathcal{B}\hat{D}_{n}}\nonumber\\
&&+\big(B_{2}(\frac{1}{4})+B_{2}(\frac{3}{4})\big)\langle e_{\lbr a^{2k}\rbr}e_{\lbr b\rbr}e_{\lbr b\rbr}\rangle_{0,3}^{\mathcal{B}\hat{D}_{n}}
+\big(B_{2}(\frac{1}{4})+B_{2}(\frac{3}{4})\big)\langle e_{\lbr a^{2k}\rbr}e_{\lbr ab\rbr}e_{\lbr ab\rbr}\rangle_{0,3}^{\mathcal{B}\hat{D}_{n}}\nonumber\\
&=& B_{2}\cdot\frac{6}{24}-B_{2}(\frac{2k}{2n-4})\cdot\frac{6}{24}+\frac{1}{2}B_{2}(\frac{n-2-k}{2n-4})\\
&&+\frac{1}{2}B_{2}(\frac{k}{2n-4})
+B_{2}(\frac{1}{4})+B_{2}(\frac{3}{4})\nonumber\\
&=&0.\nonumber
\eea
where for (\ref{31}) we have used the string equation and apply the proposition 3.4 and lemma 3.5 of \cite{JK} to obtain\footnote{Of course, we can still use Quantum Riemann-Roch  to compute this correlator.}
\ben
&&\langle e_{\lbr a^{2k}\rbr}\bar{\psi}_{1}^{1}\rangle_{1,1}^{\mathcal{B}\hat{D}_{n}}=((2n-4)\langle e_{\lbr a^{2k}\rbr}e_{\lbr a^{n-2-k}\rbr}e_{\lbr a^{n-2-k}\rbr}\rangle_{0,3}^{\mathcal{B}\hat{D}_{n}}+(2n-4)\langle e_{\lbr a^{2k}\rbr}e_{\lbr a^{k}\rbr}e_{\lbr a^{k}\rbr}\rangle_{0,3}^{\mathcal{B}\hat{D}_{n}}\\
&&+4\langle e_{\lbr a^{2k}\rbr}e_{\lbr b\rbr}e_{\lbr b\rbr}\rangle_{0,3}^{\mathcal{B}\hat{D}_{n}}
+4\langle e_{\lbr a^{2k}\rbr}e_{\lbr ab\rbr}e_{\lbr ab\rbr}\rangle_{0,3}^{\mathcal{B}\hat{D}_{n}})\langle\psi\rangle_{1,1}=\frac{6}{24}.
\een
For $k=\frac{n-2}{2}$ (thus $2|n$),
\ben
&&\langle e_{\lbr a^{n-2}\rbr}\rangle_{1,1}^{\mathcal{B}\hat{D}_{n}}\\
&=&B_{2}\langle e_{\lbr a^{n-2}\rbr}(e_{\lbr 1\rbr}\bar{\psi}_{2}^{2})\rangle_{1,2}^{\mathcal{B}\hat{D}_{n}}-B_{2}(\frac{n-2}{2n-4})\langle e_{\lbr a^{n-2}\rbr}\bar{\psi}_{1}^{1}\rangle_{1,1}^{\mathcal{B}\hat{D}_{n}}\\
&&+(n-2)B_{2}(\frac{\frac{n-2}{2}}{2n-4})\langle e_{\lbr a^{n-2}\rbr}e_{\lbr a^{\frac{n-2}{2}}\rbr}e_{\lbr a^{\frac{n-2}{2}}\rbr}\rangle_{0,3}^{\mathcal{B}\hat{D}_{n}}\\
&&+\big(B_{2}(\frac{1}{4})+B_{2}(\frac{3}{4})\big)\langle e_{\lbr a^{n-2}\rbr}e_{\lbr b\rbr}e_{\lbr b\rbr}\rangle_{0,3}^{\mathcal{B}\hat{D}_{n}}
+\big(B_{2}(\frac{1}{4})+B_{2}(\frac{3}{4})\big)\langle e_{\lbr a^{n-2}\rbr}e_{\lbr ab\rbr}e_{\lbr ab\rbr}\rangle_{0,3}^{\mathcal{B}\hat{D}_{n}}\\
&=& B_{2}\cdot\frac{3}{24}-B_{2}(\frac{n-2}{2n-4})\cdot\frac{3}{24}+\frac{B_{2}(\frac{\frac{n-2}{2}}{2n-4})}{2}
+\frac{1}{2}B_{2}(\frac{1}{4})+\frac{1}{2}B_{2}(\frac{3}{4})\\
&=&0,
\een
where we have used
\ben
&&\langle e_{\lbr a^{n-2}\rbr}\bar{\psi}_{1}^{1}\rangle_{1,1}^{\mathcal{B}\hat{D}_{n}}=((2n-4)\langle e_{\lbr a^{n-2}\rbr}e_{\lbr a^{\frac{n-2}{2}}\rbr}e_{\lbr a^{\frac{n-2}{2}}\rbr}\rangle_{0,3}^{\mathcal{B}\hat{D}_{n}}\\
&&+4\langle e_{\lbr a^{n-2}\rbr}e_{\lbr b\rbr}e_{\lbr b\rbr}\rangle_{0,3}^{\mathcal{B}\hat{D}_{n}}
+4\langle e_{\lbr a^{n-2}\rbr}e_{\lbr ab\rbr}e_{\lbr ab\rbr}\rangle_{0,3}^{\mathcal{B}\hat{D}_{n}})\langle\psi\rangle_{1,1}=
\frac{3}{24},
\een
since $2|n$.\\
Also,
\ben
&&\langle e_{\lbr b\rbr}\rangle_{1,1}^{\mathcal{B}\hat{D}_{n}}\\
&=&B_{2}\langle e_{\lbr b\rbr}(e_{\lbr 1\rbr}\bar{\psi}_{2}^{2})\rangle_{1,2}^{\mathcal{B}\hat{D}_{n}}-\frac{B_{2}(\frac{1}{4})+B_{2}(\frac{1}{4})}{2}\langle e_{\lbr b\rbr}\bar{\psi}_{1}^{1}\rangle_{1,1}^{\mathcal{B}\hat{D}_{n}}\\
&=&0,
\een
where we have use $\langle e_{\lbr b\rbr}\bar{\psi}_{1}^{1}\rangle_{1,1}^{\mathcal{B}\hat{D}_{n}}=0$ since there exist no nonzero correlators of the form $\langle e_{\lbr b\rbr}e_{\lbr x\rbr}e_{\lbr x\rbr}\rangle_{0,3}^{\mathcal{B}\hat{D}_{n}}$.  For the same reason, $\langle e_{\lbr ab\rbr}\rangle_{1,1}^{\mathcal{B}\hat{D}_{n}}$ and $\langle e_{\lbr a^{2k+1}\rbr}\rangle_{1,1}^{\mathcal{B}\hat{D}_{n}}$ vanish.

\subsection{ $\hat{E}_{6}$ and $\hat{E}_{7}$  }
The binary tetrahedral group $\hat{E}_{6}$ has order $|\hat{E}_{6}|$=24. Its generators and relations are given by
$$\hat{E}_{6}=\{a^{h}b^{j}c^{l}, 0\leq h<4, 0\leq j<2, 0\leq l<3\}/(ba=a^{-1}b, (ac)^{2}=a^{2}b, cb=a^{2}c).$$

The binary octahedral group $\hat{E}_{7}$ has order $|\hat{E}_{7}|$=48. Its generators and relations are given by
$$\hat{E}_{7}=\{a^{h}b^{j}c^{l}, 0\leq h<8, 0\leq j<2, 0\leq l<3\}/(ba=a^{-1}b, (ac)^{2}=a^{6}b, cb=a^{2}c).$$
We write $a^4=b^2=c^3=-1$. There are eight conjugacy classes $e_{\lbr 1\rbr}$, $e_{\lbr -1\rbr}$, $e_{\lbr ab\rbr}$, $e_{\lbr b\rbr}$, $e_{\lbr c^2\rbr}$, $e_{\lbr c\rbr}$, $e_{\lbr a\rbr}$, $e_{\lbr a^3\rbr}$. The nonzero genus zero length three correlators are
\ben
\langle e_{\lbr 1\rbr}e_{\lbr 1\rbr}e_{\lbr 1\rbr}\rangle_{0,3}^{\mathcal{B}\hat{E}_{7}}=\frac{1}{48},&
\langle e_{\lbr 1\rbr}e_{\lbr -1\rbr}e_{\lbr -1\rbr}\rangle_{0,3}^{\mathcal{B}\hat{E}_{7}}=\frac{1}{48},&
\langle e_{\lbr 1\rbr}e_{\lbr ab\rbr}e_{\lbr ab\rbr}\rangle_{0,3}^{\mathcal{B}\hat{E}_{7}}=\frac{1}{4},
\een
\ben
\langle e_{\lbr 1\rbr}e_{\lbr b\rbr}e_{\lbr b\rbr}\rangle_{0,3}^{\mathcal{B}\hat{E}_{7}}
=\langle e_{\lbr 1\rbr}e_{\lbr a\rbr}e_{\lbr a\rbr}\rangle_{0,3}^{\mathcal{B}\hat{E}_{7}}
=\langle e_{\lbr 1\rbr}e_{\lbr a^3\rbr}e_{\lbr a^3\rbr}\rangle_{0,3}^{\mathcal{B}\hat{E}_{7}}=\frac{1}{8},
\een
\ben
\langle e_{\lbr 1\rbr}e_{\lbr c^2\rbr}e_{\lbr c^2\rbr}\rangle_{0,3}^{\mathcal{B}\hat{E}_{7}}
=\langle e_{\lbr 1\rbr}e_{\lbr c\rbr}e_{\lbr c\rbr}\rangle_{0,3}^{\mathcal{B}\hat{E}_{7}}=\frac{1}{6},&
\langle e_{\lbr -1\rbr}e_{\lbr ab\rbr}e_{\lbr ab\rbr}\rangle_{0,3}^{\mathcal{B}\hat{E}_{7}}=\frac{1}{4},\\
\langle e_{\lbr -1\rbr}e_{\lbr b\rbr}e_{\lbr b\rbr}\rangle_{0,3}^{\mathcal{B}\hat{E}_{7}}
=\langle e_{\lbr -1\rbr}e_{\lbr a\rbr}e_{\lbr a^3\rbr}\rangle_{0,3}^{\mathcal{B}\hat{E}_{7}}=\frac{1}{8},&
\langle e_{\lbr -1\rbr}e_{\lbr c^2\rbr}e_{\lbr c\rbr}\rangle_{0,3}^{\mathcal{B}\hat{E}_{7}}=\frac{1}{6},
\een
\ben
\langle e_{\lbr ab\rbr}e_{\lbr ab\rbr}e_{\lbr b\rbr}\rangle_{0,3}^{\mathcal{B}\hat{E}_{7}}=\frac{1}{2},&
\langle e_{\lbr ab\rbr}e_{\lbr ab\rbr}e_{\lbr c\rbr}\rangle_{0,3}^{\mathcal{B}\hat{E}_{7}}
=\langle e_{\lbr ab\rbr}e_{\lbr ab\rbr}e_{\lbr c^2\rbr}\rangle_{0,3}^{\mathcal{B}\hat{E}_{7}}=1,\\
\langle e_{\lbr ab\rbr}e_{\lbr b\rbr}e_{\lbr a\rbr}\rangle_{0,3}^{\mathcal{B}\hat{E}_{7}}
=\langle e_{\lbr ab\rbr}e_{\lbr b\rbr}e_{\lbr a^3\rbr}\rangle_{0,3}^{\mathcal{B}\hat{E}_{7}}=\frac{1}{2},&
\langle e_{\lbr ab\rbr}e_{\lbr c^2\rbr}e_{\lbr a\rbr}\rangle_{0,3}^{\mathcal{B}\hat{E}_{7}}
=\langle e_{\lbr ab\rbr}e_{\lbr c^2\rbr}e_{\lbr a^3\rbr}\rangle_{0,3}^{\mathcal{B}\hat{E}_{7}}=\frac{1}{2},\\
\langle e_{\lbr ab\rbr}e_{\lbr c\rbr}e_{\lbr a\rbr}\rangle_{0,3}^{\mathcal{B}\hat{E}_{7}}
=\langle e_{\lbr ab\rbr}e_{\lbr c\rbr}e_{\lbr a^3\rbr}\rangle_{0,3}^{\mathcal{B}\hat{E}_{7}}=\frac{1}{2},&
\langle e_{\lbr b\rbr}e_{\lbr b\rbr}e_{\lbr b\rbr}\rangle_{0,3}^{\mathcal{B}\hat{E}_{7}}=\frac{1}{2},
\een
\ben
\langle e_{\lbr b\rbr}e_{\lbr c^2\rbr}e_{\lbr c^2\rbr}\rangle_{0,3}^{\mathcal{B}\hat{E}_{7}}
=\langle e_{\lbr b\rbr}e_{\lbr c^2\rbr}e_{\lbr c\rbr}\rangle_{0,3}^{\mathcal{B}\hat{E}_{7}}
=\langle e_{\lbr b\rbr}e_{\lbr c\rbr}e_{\lbr c\rbr}\rangle_{0,3}^{\mathcal{B}\hat{E}_{7}}=\frac{1}{2},\\
\langle e_{\lbr b\rbr}e_{\lbr a\rbr}e_{\lbr a\rbr}\rangle_{0,3}^{\mathcal{B}\hat{E}_{7}}
=\langle e_{\lbr b\rbr}e_{\lbr a\rbr}e_{\lbr a^3\rbr}\rangle_{0,3}^{\mathcal{B}\hat{E}_{7}}
=\langle e_{\lbr b\rbr}e_{\lbr a^3\rbr}e_{\lbr a^3\rbr}\rangle_{0,3}^{\mathcal{B}\hat{E}_{7}}=\frac{1}{8},
\een
\ben
\langle e_{\lbr c^2\rbr}e_{\lbr c^2\rbr}e_{\lbr c\rbr}\rangle_{0,3}^{\mathcal{B}\hat{E}_{7}}=\frac{1}{2},&
\langle e_{\lbr c^2\rbr}e_{\lbr a\rbr}e_{\lbr a^3\rbr}\rangle_{0,3}^{\mathcal{B}\hat{E}_{7}}=\frac{1}{2},&
\langle e_{\lbr c\rbr}e_{\lbr c\rbr}e_{\lbr c\rbr}\rangle_{0,3}^{\mathcal{B}\hat{E}_{7}}=\frac{1}{2},
\een
\ben
\langle e_{\lbr c^2\rbr}e_{\lbr c^2\rbr}e_{\lbr c^2\rbr}\rangle_{0,3}^{\mathcal{B}\hat{E}_{7}}
=\langle e_{\lbr c^2\rbr}e_{\lbr c\rbr}e_{\lbr c\rbr}\rangle_{0,3}^{\mathcal{B}\hat{E}_{7}}=\frac{1}{6},\\
\langle e_{\lbr c\rbr}e_{\lbr a\rbr}e_{\lbr a\rbr}\rangle_{0,3}^{\mathcal{B}\hat{E}_{7}}
=\langle e_{\lbr c\rbr}e_{\lbr a^3\rbr}e_{\lbr a^3\rbr}\rangle_{0,3}^{\mathcal{B}\hat{E}_{7}}=\frac{1}{2}.
\een

The quadratic operator is
\ben
&&\big(\frac{A_{p+1}(V_{\rho_{1}})z^p}{(p+1)!}\big)^{\wedge}\\
&=&\frac{2B_{p+1}}{(p+1)!}\pd_{\lbr 1\rbr,1+p}
-\frac{2B_{p+1}}{(p+1)!}\sum_{l= 0}^{\infty}t_{l}^{\lbr 1\rbr}\pd_{\lbr 1\rbr,l+p}
-\frac{2B_{p+1}(\frac{1}{2})}{(p+1)!}\sum_{l= 0}^{\infty}t_{l}^{\lbr -1\rbr}\pd_{\lbr -1\rbr,l+p}\\
&&-\frac{B_{p+1}(\frac{1}{4})+B_{p+1}(\frac{3}{4})}{(p+1)!}\sum_{l= 0}^{\infty}t_{l}^{\lbr ab\rbr}\pd_{\lbr ab\rbr,l+p}
-\frac{B_{p+1}(\frac{1}{4})+B_{p+1}(\frac{3}{4})}{(p+1)!}\sum_{l= 0}^{\infty}t_{l}^{\lbr b\rbr}\pd_{\lbr b\rbr,l+p}\\
&&-\frac{B_{p+1}(\frac{1}{3})+B_{p+1}(\frac{2}{3})}{(p+1)!}\sum_{l= 0}^{\infty}t_{l}^{\lbr c^2\rbr}\pd_{\lbr c^2\rbr,l+p}
-\frac{B_{p+1}(\frac{1}{6})+B_{p+1}(\frac{5}{6})}{(p+1)!}\sum_{l= 0}^{\infty}t_{l}^{\lbr c\rbr}\pd_{\lbr c\rbr,l+p}\\
&&-\frac{B_{p+1}(\frac{1}{8})+B_{p+1}(\frac{7}{8})}{(p+1)!}\sum_{l= 0}^{\infty}t_{l}^{\lbr a\rbr}\pd_{\lbr a\rbr,l+p}
-\frac{B_{p+1}(\frac{3}{8})+B_{p+1}(\frac{5}{8})}{(p+1)!}\sum_{l= 0}^{\infty}t_{l}^{\lbr a^3\rbr}\pd_{\lbr a^3\rbr,l+p}\\
&&+\frac{\hbar^2}{2}\sum_{l= 0}^{p-1}(-1)^l\Bigg(
48\cdot\frac{2B_{p+1}}{(p+1)!}\pd_{\lbr 1\rbr,l}\pd_{\lbr 1\rbr,p-1-l}
+48\cdot\frac{2B_{p+1}(\frac{1}{2})}{(p+1)!}\pd_{\lbr -1 \rbr,l}\pd_{\lbr -1\rbr,p-1-l}
\\
&&+4\cdot\frac{B_{p+1}(\frac{1}{4})+B_{p+1}(\frac{3}{4})}{(p+1)!}\pd_{\lbr ab\rbr,l}\pd_{\lbr ab\rbr,p-1-l}
+8\cdot\frac{B_{p+1}(\frac{1}{4})+B_{p+1}(\frac{3}{4})}{(p+1)!}\pd_{\lbr b\rbr,l}\pd_{\lbr b\rbr,p-1-l}\\
&&+6\cdot\frac{B_{p+1}(\frac{1}{3})+B_{p+1}(\frac{2}{3})}{(p+1)!}\pd_{\lbr c^2\rbr,l}\pd_{\lbr c^2\rbr,p-1-l}
+6\cdot\frac{B_{p+1}(\frac{1}{6})+B_{p+1}(\frac{5}{6})}{(p+1)!}\pd_{\lbr c\rbr,l}\pd_{\lbr c\rbr,p-1-l}\\
&&+8\cdot\frac{B_{p+1}(\frac{1}{8})+B_{p+1}(\frac{7}{8})}{(p+1)!}\pd_{\lbr a\rbr,l}\pd_{\lbr a\rbr,p-1-l}
+8\cdot\frac{B_{p+1}(\frac{3}{8})+B_{p+1}(\frac{5}{8})}{(p+1)!}\pd_{\lbr a^3\rbr,l}\pd_{\lbr a^3\rbr,p-1-l}
\Bigg).
\een
Thus
\ben
&&\langle \ch_{1}e_{\lbr b\rbr}\rangle_{1,1}^{\mathcal{B}\hat{E}_7}\\
&=& B_2\langle e_{\lbr b\rbr}(e_{\lbr 1\rbr}\bar{\psi}_{2}^{2})\rangle_{1,2}^{\mathcal{B}\hat{E}_7}
-\frac{B_{2}(\frac{1}{4})+B_{2}(\frac{3}{4})}{2}\langle e_{\lbr b\rbr}\bar{\psi}_{1}^{1})\rangle_{1,1}^{\mathcal{B}\hat{E}_7}+2\cdot\frac{B_{2}(\frac{1}{4})+B_{2}(\frac{3}{4})}{2}\langle e_{\lbr b\rbr}e_{\lbr ab\rbr}e_{\lbr ab\rbr}\rangle_{0,3}^{\mathcal{B}\hat{E}_7}\\
&&+4\cdot\frac{B_{2}(\frac{1}{4})+B_{2}(\frac{3}{4})}{2}\langle e_{\lbr b\rbr}e_{\lbr b\rbr}e_{\lbr b\rbr}\rangle_{0,3}^{\mathcal{B}\hat{E}_7}+3\cdot\frac{B_{2}(\frac{1}{3})
+B_{2}(\frac{2}{3})}{2}\langle e_{\lbr b\rbr}e_{\lbr c^{2}\rbr}e_{\lbr c^{2}\rbr}\rangle_{0,3}^{\mathcal{B}\hat{E}_7}
+3\cdot\frac{B_{2}(\frac{1}{6})
+B_{2}(\frac{5}{6})}{2}\langle e_{\lbr b\rbr}e_{\lbr c\rbr}e_{\lbr c\rbr}\rangle_{0,3}^{\mathcal{B}\hat{E}_7}\\
&&+4\cdot\frac{B_{2}(\frac{1}{8})
+B_{2}(\frac{7}{8})}{2}\langle e_{\lbr b\rbr}e_{\lbr a\rbr}e_{\lbr a\rbr}\rangle_{0,3}^{\mathcal{B}\hat{E}_7}
+4\cdot\frac{B_{2}(\frac{3}{8})
+B_{2}(\frac{5}{8})}{2}\langle e_{\lbr b\rbr}e_{\lbr a^{3}\rbr}e_{\lbr a^{3}\rbr}\rangle_{0,3}^{\mathcal{B}\hat{E}_7}\\
&=& \frac{14}{24}B_2-\frac{14}{24}B_2(\frac{1}{4})+B_{2}(\frac{1}{4})+2B_{2}(\frac{1}{4})+\frac{3}{2}B_{2}(\frac{1}{3})\\
&&+\frac{3}{2}B_{2}(\frac{1}{6})+\frac{1}{2}B_{2}(\frac{1}{8})+\frac{1}{2}B_{2}(\frac{3}{8})\\
&=&0,
\een
where we have used
\ben
\langle e_{\lbr b\rbr}\bar{\psi}\rangle_{1,1}^{\mathcal{B}\hat{E}_7}&=&\Big(4\langle e_{\lbr b\rbr}e_{\lbr ab\rbr}e_{\lbr ab\rbr}\rangle_{0,3}^{\mathcal{B}\hat{E}_7}+8\langle e_{\lbr b\rbr}e_{\lbr b\rbr}e_{\lbr b\rbr}\rangle_{0,3}^{\mathcal{B}\hat{E}_7}+6\langle e_{\lbr b\rbr}e_{\lbr c^{2}\rbr}e_{\lbr c^{2}\rbr}\rangle_{0,3}^{\mathcal{B}\hat{E}_7}\\
&&+6\langle e_{\lbr b\rbr}e_{\lbr c\rbr}e_{\lbr c\rbr}\rangle_{0,3}^{\mathcal{B}\hat{E}_7}
+8\langle e_{\lbr b\rbr}e_{\lbr a\rbr}e_{\lbr a\rbr}\rangle_{0,3}^{\mathcal{B}\hat{E}_7}
+8\langle e_{\lbr b\rbr}e_{\lbr a^{3}\rbr}e_{\lbr a^{3}\rbr}\rangle_{0,3}^{\mathcal{B}\hat{E}_7}
\Big)\cdot\frac{1}{24}
=\frac{14}{24}.
\een

\subsection{The binary icosahedral group $\hat{E}_{8}$ }
The group $\hat{E}_{8}=\{a,b|a^5=b^3=(ba)^2=-1\}$. By \cite{JK}, when $\rho$ runns over the irreducible representations of $\hat{E}_{8}$, the linear combinations
\ben
f_{\rho}=\sum_{\lbr \gamma\rbr}\frac{\chi_{\rho}(1)}{|G|}\chi_{\rho}(\lbr\gamma\rbr)e_{\lbr\gamma\rbr}
\een
form a semisimple basis for the quantum cohomology of $\mathcal{B}\hat{E}_{8}$. Thus the genus zero length three correlators are easily computed by inverting the character table (see e.g. \cite{Ste} ). We have the transition matrix
\ben
\begin{pmatrix}e_{\lbr 1\rbr}\\e_{\lbr -1\rbr}\\e_{\lbr a\rbr}\\e_{\lbr a^2\rbr}\\
e_{\lbr a^3\rbr}\\e_{\lbr a^4\rbr}\\e_{\lbr b\rbr}\\e_{\lbr b^2\rbr}\\e_{\lbr ab\rbr}\end{pmatrix}
=
\begin{pmatrix}
1& 1& 1& 1& 1& 1& 1& 1& 1\\
1& -1& -1& 1& 1& -1& 1& 1& -1\\
12 &3(1+\sqrt{5}) & 3(1-\sqrt{5}) & 2(1+\sqrt{5})&2(1-\sqrt{5})& 3& -3& 0& -2\\
12 &3(-1+\sqrt{5}) & 3(-1-\sqrt{5}) & 2(1-\sqrt{5})&2(1+\sqrt{5})& -3& -3& 0& 2\\
12 &3(1-\sqrt{5}) & 3(1+\sqrt{5}) & 2(1-\sqrt{5})&2(1+\sqrt{5})& 3& -3& 0& -2\\
12 &3(-1-\sqrt{5}) & 3(-1+\sqrt{5}) & 2(1+\sqrt{5})&2(1-\sqrt{5})& -3& -3& 0& 2\\
20 & 10& 10& 0& 0& -5& 5& -4& 0 \\
20 & -10& -10& 0& 0& 5& 5& -4& 0 \\
30& 0& 0& -10& -10& 0& 0& 6& 0
\end{pmatrix}
\begin{pmatrix}
f_{\chi_{1}}\\f_{\chi_{2}}\\f_{\chi_{3}}\\f_{\chi_{4}}\\f_{\chi_{5}}\\f_{\chi_{6}}\\f_{\chi_{7}}\\f_{\chi_{8}}\\f_{\chi_{9}}
\end{pmatrix}.
\een
Denote this transition matrix by $(C_{\gamma}^{\alpha})$, then
\ben
\langle e_{\gamma_1} e_{\gamma_2}e_{\gamma_3}\rangle^{\mathcal{B}\hat{E}_{8}}
=\sum_{i=1}^{9}C_{\gamma_1}^{\chi_i}C_{\gamma_2}^{\chi_i}C_{\gamma_3}^{\chi_i}\big(\frac{\chi_{i}(1)}{120}\big)^2.
\een
Thus the nonzero genus zero length three correlators are
\begin{flushleft}
\ben
\langle e_{\lbr 1\rbr}e_{\lbr 1\rbr}e_{\lbr 1\rbr}\rangle_{0,3}^{\mathcal{B}\hat{E}_8}=\frac{1}{120},
&\langle e_{\lbr 1\rbr}e_{\lbr -1\rbr}e_{\lbr -1\rbr}\rangle_{0,3}^{\mathcal{B}\hat{E}_8}=\frac{1}{120},\\
\langle e_{\lbr 1\rbr}e_{\lbr a^{i}\rbr}e_{\lbr a^{i}\rbr}\rangle_{0,3}^{\mathcal{B}\hat{E}_8}=\frac{1}{10}, \text{ for }
i=1,2,3,4,
&\langle e_{\lbr 1\rbr}e_{\lbr b\rbr}e_{\lbr b\rbr}\rangle_{0,3}^{\mathcal{B}\hat{E}_8}
=\langle e_{\lbr 1\rbr}e_{\lbr b^2\rbr}e_{\lbr b^2\rbr}\rangle_{0,3}^{\mathcal{B}\hat{E}_8}=\frac{1}{6},\\
\langle e_{\lbr 1\rbr}e_{\lbr ab\rbr}e_{\lbr ab\rbr}\rangle_{0,3}^{\mathcal{B}\hat{E}_8}=\frac{1}{4},
&\langle e_{\lbr -1\rbr}e_{\lbr a\rbr}e_{\lbr a^4\rbr}\rangle_{0,3}^{\mathcal{B}\hat{E}_8}
=\langle e_{\lbr -1\rbr}e_{\lbr a^2\rbr}e_{\lbr a^3\rbr}\rangle_{0,3}^{\mathcal{B}\hat{E}_8}=\frac{1}{10},\\
\langle e_{\lbr -1\rbr}e_{\lbr b\rbr}e_{\lbr b^2\rbr}\rangle_{0,3}^{\mathcal{B}\hat{E}_8}=\frac{1}{6},
&\langle e_{\lbr -1\rbr}e_{\lbr ab\rbr}e_{\lbr ab\rbr}\rangle_{0,3}^{\mathcal{B}\hat{E}_8}=\frac{1}{4},
\een
\ben
&&\langle e_{\lbr a^{i}\rbr}e_{\lbr a^{j}\rbr}e_{\lbr a^{k}\rbr}\rangle_{0,3}^{\mathcal{B}\hat{E}_8}=\frac{1}{10}, \text{ for }
(i,j,k)=(1,1,2),(1,2,3),(1,3,4),(2,2,4),(2,4,4),(3,3,4);\\
&&\langle e_{\lbr a^{i}\rbr}e_{\lbr a^{j}\rbr}e_{\lbr a^{k}\rbr}\rangle_{0,3}^{\mathcal{B}\hat{E}_8}=\frac{1}{2}, \text{ for }
(i,j,k)=(1,1,1),(1,4,4),(2,2,3),(3,3,3);\\
&&\langle e_{\lbr a^{i}\rbr}e_{\lbr a^{j}\rbr}e_{\lbr b\rbr}\rangle_{0,3}^{\mathcal{B}\hat{E}_8}=\frac{1}{2}, \text{ for }
(i,j)=(1,1),(1,2),(2,2),(3,3),(3,4),(4,4);\\
&&\langle e_{\lbr a^{i}\rbr}e_{\lbr a^{j}\rbr}e_{\lbr b^2\rbr}\rangle_{0,3}^{\mathcal{B}\hat{E}_8}=\frac{1}{2}, \text{ for }
(i,j)=(1,3),(1,4),(2,3),(2,4);\\
&&\langle e_{\lbr a^{i}\rbr}e_{\lbr b\rbr}e_{\lbr b\rbr}\rangle_{0,3}^{\mathcal{B}\hat{E}_8}=
\langle e_{\lbr a^{i}\rbr}e_{\lbr b^2\rbr}e_{\lbr b^2\rbr}\rangle_{0,3}^{\mathcal{B}\hat{E}_8}=\frac{1}{2}, \text{ for }
i=1,3;\\
&&\langle e_{\lbr a^{i}\rbr}e_{\lbr b\rbr}e_{\lbr b^2\rbr}\rangle_{0,3}^{\mathcal{B}\hat{E}_8}=\frac{1}{2}, \text{ for }
i=2,4;\\
&&\langle e_{\lbr a^{i}\rbr}e_{\lbr b\rbr}e_{\lbr ab\rbr}\rangle_{0,3}^{\mathcal{B}\hat{E}_8}=
\langle e_{\lbr a^{i}\rbr}e_{\lbr b^2\rbr}e_{\lbr ab\rbr}\rangle_{0,3}^{\mathcal{B}\hat{E}_8}=\frac{1}{2}, \text{ for }
i=1,2,3,4;\\
&&\langle e_{\lbr a^{i}\rbr}e_{\lbr ab\rbr}e_{\lbr ab\rbr}\rangle_{0,3}^{\mathcal{B}\hat{E}_8}=1, \text{ for }
i=1,2,3,4 ;\\
&&\langle e_{\lbr b\rbr}e_{\lbr b\rbr}e_{\lbr b^2\rbr}\rangle_{0,3}^{\mathcal{B}\hat{E}_8}
=\langle e_{\lbr b^2\rbr}e_{\lbr b^2\rbr}e_{\lbr b^2\rbr}\rangle_{0,3}^{\mathcal{B}\hat{E}_8}=\frac{1}{6},\\
&&\langle e_{\lbr b\rbr}e_{\lbr b\rbr}e_{\lbr b\rbr}\rangle_{0,3}^{\mathcal{B}\hat{E}_8}
=\langle e_{\lbr b\rbr}e_{\lbr b\rbr}e_{\lbr ab\rbr}\rangle_{0,3}^{\mathcal{B}\hat{E}_8}
=\langle e_{\lbr b\rbr}e_{\lbr b^2\rbr}e_{\lbr b^2\rbr}\rangle_{0,3}^{\mathcal{B}\hat{E}_8}
=\langle e_{\lbr b\rbr}e_{\lbr b^2\rbr}e_{\lbr ab\rbr}\rangle_{0,3}^{\mathcal{B}\hat{E}_8}\\
&=&\langle e_{\lbr b\rbr}e_{\lbr ab\rbr}e_{\lbr ab\rbr}\rangle_{0,3}^{\mathcal{B}\hat{E}_8}
=\langle e_{\lbr b^2\rbr}e_{\lbr b^2\rbr}e_{\lbr ab\rbr}\rangle_{0,3}^{\mathcal{B}\hat{E}_8}
=\langle e_{\lbr b^2\rbr}e_{\lbr ab\rbr}e_{\lbr ab\rbr}\rangle_{0,3}^{\mathcal{B}\hat{E}_8}
=\langle e_{\lbr ab\rbr}e_{\lbr ab\rbr}e_{\lbr ab\rbr}\rangle_{0,3}^{\mathcal{B}\hat{E}_8}=1.
\een
\end{flushleft}
Take $\chi_2$ as the standard representation of $\hat{E}_8$.
\ben
&&\big(\frac{A_{p+1}(V_{\chi_{2}})z^p}{(p+1)!}\big)^{\wedge}\\
&=&\frac{2B_{p+1}}{(p+1)!}\pd_{\lbr 1\rbr,1+p}
-\frac{2B_{p+1}}{(p+1)!}\sum_{l= 0}^{\infty}t_{l}^{\lbr 1\rbr}\pd_{\lbr 1\rbr,l+p}
-\frac{2B_{p+1}(\frac{1}{2})}{(p+1)!}\sum_{l= 0}^{\infty}t_{l}^{\lbr -1\rbr}\pd_{\lbr -1\rbr,l+p}\\
&&-\frac{B_{p+1}(\frac{1}{10})+B_{p+1}(\frac{9}{10})}{(p+1)!}\sum_{l= 0}^{\infty}t_{l}^{\lbr a\rbr}\pd_{\lbr a\rbr,l+p}
-\frac{B_{p+1}(\frac{1}{5})+B_{p+1}(\frac{4}{5})}{(p+1)!}\sum_{l= 0}^{\infty}t_{l}^{\lbr a^2\rbr}\pd_{\lbr a^2\rbr,l+p}\\
&&-\frac{B_{p+1}(\frac{3}{10})+B_{p+1}(\frac{7}{10})}{(p+1)!}\sum_{l= 0}^{\infty}t_{l}^{\lbr a^3\rbr}\pd_{\lbr a^3\rbr,l+p}
-\frac{B_{p+1}(\frac{2}{5})+B_{p+1}(\frac{3}{5})}{(p+1)!}\sum_{l= 0}^{\infty}t_{l}^{\lbr a^4\rbr}\pd_{\lbr a^4\rbr,l+p}\\
&&-\frac{B_{p+1}(\frac{1}{6})+B_{p+1}(\frac{5}{6})}{(p+1)!}\sum_{l= 0}^{\infty}t_{l}^{\lbr b\rbr}\pd_{\lbr b\rbr,l+p}
-\frac{B_{p+1}(\frac{1}{3})+B_{p+1}(\frac{2}{3})}{(p+1)!}\sum_{l= 0}^{\infty}t_{l}^{\lbr b^2\rbr}\pd_{\lbr b^2\rbr,l+p}\\
&&-\frac{B_{p+1}(\frac{1}{4})+B_{p+1}(\frac{3}{4})}{(p+1)!}\sum_{l= 0}^{\infty}t_{l}^{\lbr ab\rbr}\pd_{\lbr ab\rbr,l+p}\\
&&+\frac{\hbar^2}{2}\sum_{l= 0}^{p-1}(-1)^l\Bigg(
120\cdot\frac{2B_{p+1}}{(p+1)!}\pd_{\lbr 1\rbr,l}\pd_{\lbr 1\rbr,p-1-l}
+120\cdot\frac{2B_{p+1}(\frac{1}{2})}{(p+1)!}\pd_{\lbr -1 \rbr,l}\pd_{\lbr -1\rbr,p-1-l}
\\
&&+10\cdot\frac{B_{p+1}(\frac{1}{10})+B_{p+1}(\frac{9}{10})}{(p+1)!}\pd_{\lbr a\rbr,l}\pd_{\lbr a\rbr,p-1-l}
+10\cdot\frac{B_{p+1}(\frac{1}{5})+B_{p+1}(\frac{4}{5})}{(p+1)!}\pd_{\lbr a^2\rbr,l}\pd_{\lbr a^2\rbr,p-1-l}\\
&&+10\cdot\frac{B_{p+1}(\frac{3}{10})+B_{p+1}(\frac{7}{10})}{(p+1)!}\pd_{\lbr a^3\rbr,l}\pd_{\lbr a^3\rbr,p-1-l}
+10\cdot\frac{B_{p+1}(\frac{2}{5})+B_{p+1}(\frac{3}{5})}{(p+1)!}\pd_{\lbr a^4\rbr,l}\pd_{\lbr a^4\rbr,p-1-l}\\
&&+6\cdot\frac{B_{p+1}(\frac{1}{6})+B_{p+1}(\frac{5}{6})}{(p+1)!}\pd_{\lbr b\rbr,l}\pd_{\lbr b\rbr,p-1-l}
+6\cdot\frac{B_{p+1}(\frac{1}{3})+B_{p+1}(\frac{2}{3})}{(p+1)!}\pd_{\lbr b^2\rbr,l}\pd_{\lbr b^2\rbr,p-1-l}\\
&&+4\cdot\frac{B_{p+1}(\frac{1}{4})+B_{p+1}(\frac{3}{4})}{(p+1)!}\pd_{\lbr ab\rbr,l}\pd_{\lbr ab\rbr,p-1-l}
\Bigg).
\een

\bea
&&\langle \ch_{1}^{\chi_2}e_{\lbr a\rbr}\rangle_{1,1}^{\mathcal{B}\hat{E}_8}\nonumber\\
&=& B_2\langle e_{\lbr a\rbr}(e_{\lbr 1\rbr}\bar{\psi}_{2}^{2})\rangle_{1,2}^{\mathcal{B}\hat{E}_8}
-\frac{B_{2}(\frac{1}{10})+B_{2}(\frac{9}{10})}{2}\langle e_{\lbr a\rbr}\bar{\psi}_{1}^{1})\rangle_{1,1}^{\mathcal{B}\hat{E}_8}+5\cdot\frac{B_{2}(\frac{1}{10})+B_{2}(\frac{9}{10})}{2}\langle e_{\lbr a\rbr}e_{\lbr a\rbr}e_{\lbr a\rbr}\rangle_{0,3}^{\mathcal{B}\hat{E}_8}\nonumber\\
&&+5\cdot\frac{B_{2}(\frac{2}{5})+B_{2}(\frac{3}{5})}{2}\langle e_{\lbr a\rbr}e_{\lbr a^4\rbr}e_{\lbr a^4\rbr}\rangle_{0,3}^{\mathcal{B}\hat{E}_8}+3\cdot\frac{B_{2}(\frac{1}{6})
+B_{2}(\frac{5}{6})}{2}\langle e_{\lbr a\rbr}e_{\lbr b\rbr}e_{\lbr b\rbr}\rangle_{0,3}^{\mathcal{B}\hat{E}_8}\nonumber\\
&&+3\cdot\frac{B_{2}(\frac{1}{3})
+B_{2}(\frac{2}{3})}{2}\langle e_{\lbr a\rbr}e_{\lbr b^2\rbr}e_{\lbr b^2\rbr}\rangle_{0,3}^{\mathcal{B}\hat{E}_8}+2\cdot\frac{B_{2}(\frac{1}{4})
+B_{2}(\frac{3}{4})}{2}\langle e_{\lbr a\rbr}e_{\lbr ab\rbr}e_{\lbr ab\rbr}\rangle_{0,3}^{\mathcal{B}\hat{E}_8}\nonumber\\
&=& \frac{20}{24}B_2 -\frac{20}{24}B_2(\frac{1}{10})+5B_{2}(\frac{1}{10})\cdot\frac{1}{2}\nonumber\\
&&+5B_{2}(\frac{2}{5})\cdot\frac{1}{2}+3B_2(\frac{1}{6})\cdot\frac{1}{2}+3B_2(\frac{1}{3})\cdot\frac{1}{2}+2B_2(\frac{1}{4})\\
&=&0,\nonumber
\eea
where we have used
\ben
\langle e_{\lbr a\rbr}\bar{\psi}\rangle_{1,1}^{\mathcal{B}\hat{E}_8}&=&\Big(10\langle e_{\lbr a\rbr}e_{\lbr a\rbr}e_{\lbr a\rbr}\rangle_{0,3}^{\mathcal{B}\hat{E}_8}+10\langle e_{\lbr a\rbr}e_{\lbr a^{4}\rbr}e_{\lbr a^{4}\rbr}\rangle_{0,3}^{\mathcal{B}\hat{E}_8}\\
&&+6\langle e_{\lbr a\rbr}e_{\lbr b\rbr}e_{\lbr b\rbr}\rangle_{0,3}^{\mathcal{B}\hat{E}_8}+6\langle e_{\lbr a\rbr}e_{\lbr b^{2}\rbr}e_{\lbr b^{2}\rbr}\rangle_{0,3}^{\mathcal{B}\hat{E}_8}+4\langle e_{\lbr a\rbr}e_{\lbr ab\rbr}e_{\lbr ab\rbr}\rangle_{0,3}^{\mathcal{B}\hat{E}_8}\Big)\langle\psi\rangle_{1,1}\\
&=&(5+5+3+3+4)\cdot\frac{1}{24}=\frac{20}{24}.
\een

\ben
&&\langle \ch_{1}^{\chi_2}e_{\lbr a^{2}\rbr}\rangle_{1,1}^{\mathcal{B}\hat{E}_8}\\
&=& B_2\langle e_{\lbr a^{2}\rbr}(e_{\lbr 1\rbr}\bar{\psi}_{2}^{2})\rangle_{1,2}^{\mathcal{B}\hat{E}_8}
-\frac{B_{2}(\frac{1}{5})+B_{2}(\frac{4}{5})}{2}\langle e_{\lbr a^{2}\rbr}\bar{\psi}_{1}^{1})\rangle_{1,1}^{\mathcal{B}\hat{E}_8}+5\cdot\frac{B_{2}(\frac{1}{10})+B_{2}(\frac{9}{10})}{2}\langle e_{\lbr a^{2}\rbr}e_{\lbr a\rbr}e_{\lbr a\rbr}\rangle_{0,3}^{\mathcal{B}\hat{E}_8}\\
&&+5\cdot\frac{B_{2}(\frac{2}{5})+B_{2}(\frac{3}{5})}{2}\langle e_{\lbr a^{2}\rbr}e_{\lbr a^4\rbr}e_{\lbr a^4\rbr}\rangle_{0,3}^{\mathcal{B}\hat{E}_8}+2\cdot\frac{B_{2}(\frac{1}{4})
+B_{2}(\frac{3}{4})}{2}\langle e_{\lbr a\rbr}e_{\lbr ab\rbr}e_{\lbr ab\rbr}\rangle_{0,3}^{\mathcal{B}\hat{E}_8}\\
&=& \frac{6}{24}B_2-\frac{6}{24}B_2(\frac{1}{5})+5B_{2}(\frac{1}{10})\cdot\frac{1}{10}\\
&&+5B_{2}(\frac{2}{5})\cdot\frac{1}{10}+2B_2(\frac{1}{4})\\
&=&0,
\een
where we have used
\ben
\langle e_{\lbr a^{2}\rbr}\bar{\psi}\rangle_{1,1}^{\mathcal{B}\hat{E}_8}=\Big(10\langle e_{\lbr a^{2}\rbr}e_{\lbr a\rbr}e_{\lbr a\rbr}\rangle_{0,3}^{\mathcal{B}\hat{E}_8}+10\langle e_{\lbr a^{2}\rbr}e_{\lbr a^{4}\rbr}e_{\lbr a^{4}\rbr}\rangle_{0,3}^{\mathcal{B}\hat{E}_8}+4\langle e_{\lbr a^{2}\rbr}e_{\lbr ab\rbr}e_{\lbr ab\rbr}\rangle_{0,3}^{\mathcal{B}\hat{E}_8}\Big)\cdot\frac{1}{24}
=\frac{6}{24}.
\een

\ben
&&\langle \ch_{1}^{\chi_2}e_{\lbr a^{3}\rbr}\rangle_{1,1}^{\mathcal{B}\hat{E}_8}\nonumber\\
&=& B_2\langle e_{\lbr a^{3}\rbr}(e_{\lbr 1\rbr}\bar{\psi}_{2}^{2})\rangle_{1,2}^{\mathcal{B}\hat{E}_8}
-\frac{B_{2}(\frac{3}{10})+B_{2}(\frac{7}{10})}{2}\langle e_{\lbr a^{3}\rbr}\bar{\psi}_{1}^{1})\rangle_{1,1}^{\mathcal{B}\hat{E}_8}+5\cdot\frac{B_{2}(\frac{1}{5})+B_{2}(\frac{4}{5})}{2}\langle e_{\lbr a^{3}\rbr}e_{\lbr a^{2}\rbr}e_{\lbr a^{2}\rbr}\rangle_{0,3}^{\mathcal{B}\hat{E}_8}\nonumber\\
&&+5\cdot\frac{B_{2}(\frac{3}{10})+B_{2}(\frac{7}{10})}{2}\langle e_{\lbr a^{3}\rbr}e_{\lbr a^3\rbr}e_{\lbr a^3\rbr}\rangle_{0,3}^{\mathcal{B}\hat{E}_8}+3\cdot\frac{B_{2}(\frac{1}{6})
+B_{2}(\frac{5}{6})}{2}\langle e_{\lbr a^3\rbr}e_{\lbr b\rbr}e_{\lbr b\rbr}\rangle_{0,3}^{\mathcal{B}\hat{E}_8}\nonumber\\
&&+3\cdot\frac{B_{2}(\frac{1}{3})
+B_{2}(\frac{2}{3})}{2}\langle e_{\lbr a^3\rbr}e_{\lbr b^2\rbr}e_{\lbr b^2\rbr}\rangle_{0,3}^{\mathcal{B}\hat{E}_8}+2\cdot\frac{B_{2}(\frac{1}{4})
+B_{2}(\frac{3}{4})}{2}\langle e_{\lbr a^3\rbr}e_{\lbr ab\rbr}e_{\lbr ab\rbr}\rangle_{0,3}^{\mathcal{B}\hat{E}_8}\nonumber\\
&=& \frac{20}{24}B_2-\frac{20}{24}B_2(\frac{3}{10})+5B_{2}(\frac{1}{5})\cdot\frac{1}{2}\nonumber\\
&&+5B_{2}(\frac{3}{10})\cdot\frac{1}{2}+3B_2(\frac{1}{6})\cdot\frac{1}{2}+3B_2(\frac{1}{3})\cdot\frac{1}{2}+2B_2(\frac{1}{4})\\
&=&0,
\een
where we have used
\ben
\langle e_{\lbr a^3\rbr}\bar{\psi}\rangle_{1,1}^{\mathcal{B}\hat{E}_8}&=&\Big(10\langle e_{\lbr a^3\rbr}e_{\lbr a^2\rbr}e_{\lbr a^2\rbr}\rangle_{0,3}^{\mathcal{B}\hat{E}_8}+10\langle e_{\lbr a^3\rbr}e_{\lbr a^{3}\rbr}e_{\lbr a^{3}\rbr}\rangle_{0,3}^{\mathcal{B}\hat{E}_8}
+6\langle e_{\lbr a^3\rbr}e_{\lbr b\rbr}e_{\lbr b\rbr}\rangle_{0,3}^{\mathcal{B}\hat{E}_8}\\
&&+6\langle e_{\lbr a^3\rbr}e_{\lbr b^{2}\rbr}e_{\lbr b^{2}\rbr}\rangle_{0,3}^{\mathcal{B}\hat{E}_8}+4\langle e_{\lbr a^3\rbr}e_{\lbr ab\rbr}e_{\lbr ab\rbr}\rangle_{0,3}^{\mathcal{B}\hat{E}_8}\Big)\cdot\frac{1}{24}=\frac{20}{24}.
\een

\ben
&&\langle \ch_{1}^{\chi_2}e_{\lbr a^{4}\rbr}\rangle_{1,1}^{\mathcal{B}\hat{E}_8}\\
&=& B_2\langle e_{\lbr a^{4}\rbr}(e_{\lbr 1\rbr}\bar{\psi}_{2}^{2})\rangle_{1,2}^{\mathcal{B}\hat{E}_8}
-\frac{B_{2}(\frac{2}{5})+B_{2}(\frac{3}{5})}{2}\langle e_{\lbr a^{4}\rbr}\bar{\psi}_{1}^{1})\rangle_{1,1}^{\mathcal{B}\hat{E}_8}+5\cdot\frac{B_{2}(\frac{1}{5})+B_{2}(\frac{4}{5})}{2}\langle e_{\lbr a^{4}\rbr}e_{\lbr a^2\rbr}e_{\lbr a^2\rbr}\rangle_{0,3}^{\mathcal{B}\hat{E}_8}\\
&&+5\cdot\frac{B_{2}(\frac{3}{10})+B_{2}(\frac{7}{10})}{2}\langle e_{\lbr a^{4}\rbr}e_{\lbr a^3\rbr}e_{\lbr a^3\rbr}\rangle_{0,3}^{\mathcal{B}\hat{E}_8}+2\cdot\frac{B_{2}(\frac{1}{4})
+B_{2}(\frac{3}{4})}{2}\langle e_{\lbr a^4\rbr}e_{\lbr ab\rbr}e_{\lbr ab\rbr}\rangle_{0,3}^{\mathcal{B}\hat{E}_8}\\
&=& \frac{6}{24}B_2-\frac{6}{24}B_2(\frac{2}{5})+5B_{2}(\frac{1}{5})\cdot\frac{1}{10}\\
&&+5B_{2}(\frac{3}{10})\cdot\frac{1}{10}+2B_2(\frac{1}{4})\\
&=&0,
\een
where we have used
\ben
\langle e_{\lbr a^{4}\rbr}\bar{\psi}\rangle_{1,1}^{\mathcal{B}\hat{E}_8}=\Big(10\langle e_{\lbr a^{4}\rbr}e_{\lbr a^2\rbr}e_{\lbr a^2\rbr}\rangle_{0,3}^{\mathcal{B}\hat{E}_8}+10\langle e_{\lbr a^{4}\rbr}e_{\lbr a^{3}\rbr}e_{\lbr a^{3}\rbr}\rangle_{0,3}^{\mathcal{B}\hat{E}_8}+4\langle e_{\lbr a^{4}\rbr}e_{\lbr ab\rbr}e_{\lbr ab\rbr}\rangle_{0,3}^{\mathcal{B}\hat{E}_8}\Big)\cdot\frac{1}{24}
=\frac{6}{24}.
\een

\ben
&&\langle \ch_{1}^{\chi_2}e_{\lbr b\rbr}\rangle_{1,1}^{\mathcal{B}\hat{E}_8}\nonumber\\
&=& B_2\langle e_{\lbr b\rbr}(e_{\lbr 1\rbr}\bar{\psi}_{2}^{2})\rangle_{1,2}^{\mathcal{B}\hat{E}_8}
-\frac{B_{2}(\frac{1}{6})+B_{2}(\frac{5}{6})}{2}\langle e_{\lbr b\rbr}\bar{\psi}_{1}^{1})\rangle_{1,1}^{\mathcal{B}\hat{E}_8}+5\cdot\frac{B_{2}(\frac{1}{10})+B_{2}(\frac{9}{10})}{2}\langle e_{\lbr b\rbr}e_{\lbr a\rbr}e_{\lbr a\rbr}\rangle_{0,3}^{\mathcal{B}\hat{E}_8}\nonumber\\
&&+5\cdot\frac{B_{2}(\frac{1}{5})+B_{2}(\frac{4}{5})}{2}\langle e_{\lbr b\rbr}e_{\lbr a^{2}\rbr}e_{\lbr a^{2}\rbr}\rangle_{0,3}^{\mathcal{B}\hat{E}_8}\nonumber+5\cdot\frac{B_{2}(\frac{3}{10})+B_{2}(\frac{7}{10})}{2}\langle e_{\lbr b\rbr}e_{\lbr a^{3}\rbr}e_{\lbr a^{3}\rbr}\rangle_{0,3}^{\mathcal{B}\hat{E}_8}\nonumber\\
&&+5\cdot\frac{B_{2}(\frac{2}{5})+B_{2}(\frac{3}{5})}{2}\langle e_{\lbr b\rbr}e_{\lbr a^4\rbr}e_{\lbr a^4\rbr}\rangle_{0,3}^{\mathcal{B}\hat{E}_8}
+3\cdot\frac{B_{2}(\frac{1}{6})
+B_{2}(\frac{5}{6})}{2}\langle e_{\lbr b\rbr}e_{\lbr b\rbr}e_{\lbr b\rbr}\rangle_{0,3}^{\mathcal{B}\hat{E}_8}\\
&&+3\cdot\frac{B_{2}(\frac{1}{3})
+B_{2}(\frac{2}{3})}{2}\langle e_{\lbr b\rbr}e_{\lbr b^2\rbr}e_{\lbr b^2\rbr}\rangle_{0,3}^{\mathcal{B}\hat{E}_8}+2\cdot\frac{B_{2}(\frac{1}{4})
+B_{2}(\frac{3}{4})}{2}\langle e_{\lbr b\rbr}e_{\lbr ab\rbr}e_{\lbr ab\rbr}\rangle_{0,3}^{\mathcal{B}\hat{E}_8}\nonumber\\
&=& \frac{36}{24}B_2 -\frac{36}{24}B_2(\frac{1}{6})+\frac{5}{2}\Big(B_{2}(\frac{1}{10})+B_{2}(\frac{1}{5})
+B_{2}(\frac{3}{10})+B_{2}(\frac{2}{5})\Big)\\
&&+3B_2(\frac{1}{6})+3B_2(\frac{1}{3})+2B_2(\frac{1}{4})\\
&=&0,
\een
where we have used
\ben
\langle e_{\lbr b\rbr}\bar{\psi}\rangle_{1,1}^{\mathcal{B}\hat{E}_8}&=&\Big(\sum_{i=1}^{4}10\langle e_{\lbr b\rbr}e_{\lbr a^{i}\rbr}e_{\lbr a^{i}\rbr}\rangle_{0,3}^{\mathcal{B}\hat{E}_8}
+6\langle e_{\lbr b\rbr}e_{\lbr b\rbr}e_{\lbr b\rbr}\rangle_{0,3}^{\mathcal{B}\hat{E}_8}\\
&&+6\langle e_{\lbr b\rbr}e_{\lbr b^{2}\rbr}e_{\lbr b^{2}\rbr}\rangle_{0,3}^{\mathcal{B}\hat{E}_8}+4\langle e_{\lbr b\rbr}e_{\lbr ab\rbr}e_{\lbr ab\rbr}\rangle_{0,3}^{\mathcal{B}\hat{E}_8}\Big)\cdot\frac{1}{24}=\frac{36}{24}.
\een

\ben
&&\langle \ch_{1}^{\chi_2}e_{\lbr b^2\rbr}\rangle_{1,1}^{\mathcal{B}\hat{E}_8}\nonumber\\
&=& B_2\langle e_{\lbr b^2\rbr}(e_{\lbr 1\rbr}\bar{\psi}_{2}^{2})\rangle_{1,2}^{\mathcal{B}\hat{E}_8}
-\frac{B_{2}(\frac{1}{3})+B_{2}(\frac{2}{3})}{2}\langle e_{\lbr b^2\rbr}\bar{\psi}_{1}^{1})\rangle_{1,1}^{\mathcal{B}\hat{E}_8}
+3\cdot\frac{B_{2}(\frac{1}{6})
+B_{2}(\frac{5}{6})}{2}\langle e_{\lbr b^2\rbr}e_{\lbr b\rbr}e_{\lbr b\rbr}\rangle_{0,3}^{\mathcal{B}\hat{E}_8}\\
&&+3\cdot\frac{B_{2}(\frac{1}{3})
+B_{2}(\frac{2}{3})}{2}\langle e_{\lbr b^2\rbr}e_{\lbr b^2\rbr}e_{\lbr b^2\rbr}\rangle_{0,3}^{\mathcal{B}\hat{E}_8}+2\cdot\frac{B_{2}(\frac{1}{4})
+B_{2}(\frac{3}{4})}{2}\langle e_{\lbr b^2\rbr}e_{\lbr ab\rbr}e_{\lbr ab\rbr}\rangle_{0,3}^{\mathcal{B}\hat{E}_8}\nonumber\\
&=& \frac{6}{24}B_2 -\frac{6}{24}B_2(\frac{1}{3})+\frac{1}{2}B_2(\frac{1}{6})+\frac{1}{2}B_2(\frac{1}{3})+2B_2(\frac{1}{4})\\
&=&0,
\een
where we have used
\ben
\langle e_{\lbr b^2\rbr}\bar{\psi}\rangle_{1,1}^{\mathcal{B}\hat{E}_8}&=&\Big(6\langle e_{\lbr b^2\rbr}e_{\lbr b\rbr}e_{\lbr b\rbr}\rangle_{0,3}^{\mathcal{B}\hat{E}_8}+6\langle e_{\lbr b^2\rbr}e_{\lbr b^{2}\rbr}e_{\lbr b^{2}\rbr}\rangle_{0,3}^{\mathcal{B}\hat{E}_8}+4\langle e_{\lbr b^2\rbr}e_{\lbr ab\rbr}e_{\lbr ab\rbr}\rangle_{0,3}^{\mathcal{B}\hat{E}_8}\Big)\cdot\frac{1}{24}=\frac{6}{24}.
\een

\ben
&&\langle \ch_{1}^{\chi_2}e_{\lbr ab\rbr}\rangle_{1,1}^{\mathcal{B}\hat{E}_8}\nonumber\\
&=& B_2\langle e_{\lbr ab\rbr}(e_{\lbr 1\rbr}\bar{\psi}_{2}^{2})\rangle_{1,2}^{\mathcal{B}\hat{E}_8}
-\frac{B_{2}(\frac{1}{4})+B_{2}(\frac{3}{4})}{2}\langle e_{\lbr ab\rbr}\bar{\psi}_{1}^{1})\rangle_{1,1}^{\mathcal{B}\hat{E}_8}
+3\cdot\frac{B_{2}(\frac{1}{6})
+B_{2}(\frac{5}{6})}{2}\langle e_{\lbr ab\rbr}e_{\lbr b\rbr}e_{\lbr b\rbr}\rangle_{0,3}^{\mathcal{B}\hat{E}_8}\\
&&+3\cdot\frac{B_{2}(\frac{1}{3})
+B_{2}(\frac{2}{3})}{2}\langle e_{\lbr ab\rbr}e_{\lbr b^2\rbr}e_{\lbr b^2\rbr}\rangle_{0,3}^{\mathcal{B}\hat{E}_8}+2\cdot\frac{B_{2}(\frac{1}{4})
+B_{2}(\frac{3}{4})}{2}\langle e_{\lbr ab\rbr}e_{\lbr ab\rbr}e_{\lbr ab\rbr}\rangle_{0,3}^{\mathcal{B}\hat{E}_8}\nonumber\\
&=& \frac{16}{24}B_2 -\frac{16}{24}B_2(\frac{1}{4})+3B_2(\frac{1}{6})+3B_2(\frac{1}{3})+2B_2(\frac{1}{4})\\
&=&0,
\een
where we have used
\ben
\langle e_{\lbr ab\rbr}\bar{\psi}\rangle_{1,1}^{\mathcal{B}\hat{E}_8}&=&\Big(6\langle e_{\lbr ab\rbr}e_{\lbr b\rbr}e_{\lbr b\rbr}\rangle_{0,3}^{\mathcal{B}\hat{E}_8}+6\langle e_{\lbr ab\rbr}e_{\lbr b^{2}\rbr}e_{\lbr b^{2}\rbr}\rangle_{0,3}^{\mathcal{B}\hat{E}_8}+4\langle e_{\lbr ab\rbr}e_{\lbr ab\rbr}e_{\lbr ab\rbr}\rangle_{0,3}^{\mathcal{B}\hat{E}_8}\Big)\cdot\frac{1}{24}=\frac{16}{24}.
\een

\ben
&&\langle \ch_{1}^{\chi_2}e_{\lbr -1\rbr}\rangle_{1,1}^{\mathcal{B}\hat{E}_8}\\
&=& B_2\langle e_{\lbr -1\rbr}(e_{\lbr 1\rbr}\bar{\psi}_{2}^{2})\rangle_{1,2}^{\mathcal{B}\hat{E}_8}
-B_{2}(\frac{1}{2})\langle e_{\lbr -1\rbr}\bar{\psi}_{1}^{1})\rangle_{1,1}^{\mathcal{B}\hat{E}_8}+2\cdot\frac{B_{2}(\frac{1}{4})
+B_{2}(\frac{3}{4})}{2}\langle e_{\lbr -1\rbr}e_{\lbr ab\rbr}e_{\lbr ab\rbr}\rangle_{0,3}^{\mathcal{B}\hat{E}_8}\\
&=& \frac{1}{24}B_2-\frac{1}{24}B_2(\frac{1}{2})+\frac{1}{2}B_2(\frac{1}{4})\\
&=&0,
\een
where we have used
\ben
\langle e_{\lbr -1\rbr}\bar{\psi}\rangle_{1,1}^{\mathcal{B}\hat{E}_8}=4\langle e_{\lbr -1\rbr}e_{\lbr ab\rbr}e_{\lbr ab\rbr}\rangle_{0,3}^{\mathcal{B}\hat{E}_8}\cdot\frac{1}{24}
=\frac{1}{24}.
\een

\hfill\qedsymbol

\section{}
\subsection{The generalized divisor equation}
In this section we give a divisor equation for the generalized correlators, i.e. the correlators with ancestors and descendants mixed. We adopt the notations in \cite{KM} for the \emph{generalized correlators} (see also \cite{Manin}). We denote the first Chern class of the cotangent line bundle at the $i$-th marked point on $\overline{\cM}_{g,n}$ by $\phi_{i}$, and use the same notation to denote the class by pulling back $\phi_{i}$ through the absolute stabilization $st: \overline{\cM}_{g,n}(V,\beta)\rightarrow \overline{\cM}_{g,n}$ (defined for $2g-3+n\geq 0$) when no confusion should arise. The formulae are stated for the non-equivariant theory, while the statements and the proofs extends without difficulties to the equivariant theory (see the proof of theorem \ref{21}).

\begin{theorem}
(\emph{A divisor equation for generalized correlators}) Suppose $\gamma\in H^{2}(V)$. Then for $2g-2+m> 0$ we have
\bea\label{divisor}
&&\langle \gamma,\tau_{d_{1},e_{1}}\gamma_{1},\cdots,\tau_{d_{m},e_{m}}\gamma_{m}\rangle_{g,m+1,\beta}\nonumber\\
&=&(\gamma \cap \beta)\langle \tau_{d_{1},e_{1}}\gamma_{1},\cdots,\tau_{d_{m},e_{m}}\gamma_{m}\rangle_{g,m,\beta}
+\sum_{k=1}^{m}\langle \tau_{d_{1},e_{1}}\gamma_{1},\cdots,\tau_{d_{k}-1,e_{k}}(\gamma\cup\gamma_{k}),\cdots,\tau_{d_{m},e_{m}}\rangle_{g,m+1,\beta}\nonumber\\
&&+\sum_{k=1}^{m}\sum_{\beta_{1}+\beta_{2}=\beta}\sum_{a}\pm\langle \gamma, \tau_{d_{k},0}\gamma_{k},\Delta^{a}\rangle_{0,3,\beta_{1}}
\langle \tau_{d_{1},e_{1}}\gamma_{1},\cdots,\tau_{0,e_{k}-1}\Delta_{a},\cdots,\tau_{d_{m},e_{m}}\rangle_{g,m,\beta_{2}}.
\eea
Here $(\Delta_{a})$ and $(\Delta^{a})$ are poincar\'{e} dual bases of $H^{*}(V)$, and the sign arises from  permuting $\gamma_{k}$ with
$\gamma_{j}$ for $j< k$.
\end{theorem}
\emph{Proof}: For simplicity, we assume all classes are even classes. Consider the commutative diagram
\ben
\xymatrix{
  \overline{\cM}_{g,S\cup\{0\}}(V,\beta) \ar[d]_{st_{1}} \ar[r]^{f_{V}}
                &  \overline{\cM}_{g,S}(V,\beta)\ar[d]_{st_{2}}   \\
  \overline{\cM}_{g,S\cup\{0\}}\ar[r]^{f}
                &      \overline{\cM}_{g,S}
                              }
\een
in which $S=\{1,\cdots,m\}$, $f_{V}$ and $f$ are forgetting the $0$-th marked point, $st_{1}$ and $st_{2}$ are the absolute stabilization. For $j\in S$, let $D_{j}^{V}$ (resp. $D_{j}$)
be the divisor on $\overline{\cM}_{g,S\cup\{0\}}(V,\beta)$ (resp. $\overline{\cM}_{g,S\cup\{0\}}$) representing the $j$-th section of $f_{V}$ (resp. $f$). We have
\ben
\psi_{j}^{d}=(f_{V}^{*}\psi_{j})^{d}+[D_{j}^{V}]\cdot(f_{V}^{*}\psi_{j})^{d-1} ,
\een
and
\ben
\phi_{j}^{d}=(f^{*}\phi_{j})^{d}+[D_{j}]\cdot(f^{*}\phi_{j})^{d-1}
\een
for $\forall d\geq 1$ and $j\in S$. Thus
\ben
&&J_{g,S\cup\{0\}}(V,\beta)\cap ev_{0}^{*}(\gamma)ev_{S}^{*}(\alpha)\prod_{j\in S}\psi_{j}^{d_{j}}\phi_{j}^{e_{j}}
\\
&=& J_{g,S\cup\{0\}}(V,\beta)\cap ev_{0}^{*}(\gamma)ev_{S}^{*}(\alpha)\prod_{j\in S}\Big((f_{V}^{*}\psi_{j})^{d_{j}}+[D_{j}^{V}]\cdot(f_{V}^{*}\psi_{j})^{d_{j}-1}\Big)
st_{1}^{*}\Big((f^{*}\phi_{j})^{e_{j}}+[D_{j}]\cdot(f^{*}\phi_{j})^{e_{j}-1}\Big)
\\
&=& J_{g,S\cup\{0\}}(V,\beta)\cap ev_{0}^{*}(\gamma)ev_{S}^{*}(\alpha)\prod_{j\in S}(f_{V}^{*}\psi_{j})^{d_{j}}
st_{1}^{*}(f^{*}\phi_{j})^{e_{j}}\\
&&+J_{g,S\cup\{0\}}(V,\beta)\cap ev_{0}^{*}(\gamma)ev_{S}^{*}(\alpha)\sum_{k\in S, d_{k}\geq 1}[D_{k}^{V}]\prod_{j\in S}f_{V}^{*}\psi_{j}^{d_{j}-\delta_{kj}}
st_{1}^{*}(f^{*}\phi_{j})^{e_{j}}
\\
&&+J_{g,S\cup\{0\}}(V,\beta)\cap ev_{0}^{*}(\gamma)ev_{S}^{*}(\alpha)\sum_{k\in S, d_{k}\geq 1}st_{1}^{*}[D_{k}]\prod_{j\in S}\Big((f_{V}^{*}\psi_{j})^{d_{j}}+[D_{j}^{V}]\cdot(f_{V}^{*}\psi_{j})^{d_{j}-1}\Big)
st_{1}^{*}(f^{*}\phi_{j})^{e_{j}-\delta_{kj}}\\
&=&J_{g,S\cup\{0\}}(V,\beta)\cap ev_{0}^{*}(\gamma)(f_{V})^{*}\Big(\alpha\prod_{j\in S}\psi_{j}^{d_{j}}
\phi_{j}^{e_{j}}\Big)
+J_{g,S\cup\{0\}}(V,\beta)\cap ev_{0}^{*}(\gamma)\sum_{k\in S, d_{k}\geq 1}[D_{k}^{V}](f_{V})^{*}\Big(\alpha\prod_{j\in S}\psi_{j}^{d_{j}-\delta_{kj}}\phi_{j}^{e_{j}}\Big)\\
&&+J_{g,S\cup\{0\}}(V,\beta)\cap ev_{0}^{*}(\gamma)ev_{S}^{*}(\alpha)\sum_{k\in S, d_{k}\geq 1}st_{1}^{*}[D_{k}]\prod_{j\in S}\psi_{j}^{d_{j}}
\phi_{j}^{e_{j}-\delta_{kj}}.\\
\een
We push forward the first and the second summands in the last expression by $f_{V}$. Since $(f_{V})^{*}J_{g,S}(V,\beta)=J_{g,S\cup\{0\}}(V,\beta)$, the projection formula gives the first and the second terms in (\ref{divisor}). For the third summand, we apply the Proposition 1.2.1 in \cite{Manin} (See also \cite{BF}) and obtain the third term in (\ref{divisor}).\hfill\qedsymbol\\

Once one has known (\ref{divisor}), one can give\\

\emph{Another proof}: For simplicity we assume every $(d_{i},e_{i})=(0,0)$ except for one. We also assume all $\gamma_{i}$ are even classes, so that we can ignore the signs. The general cases go the same way. Thus we need to show
\bea\label{16}
&&\langle \gamma,\gamma_{1},\cdots,\gamma_{m},\gamma_{m+1}\psi^{d}\phi^{e}\rangle_{g,m+2,\beta}\nonumber\\
&=& (\gamma \cap \beta)\langle \gamma_{1},\cdots,\gamma_{m},\gamma_{m+1}\psi^{d}\phi^{e}\rangle_{g,m+1,\beta}
+\langle \gamma\cup\gamma_{m+1}\psi^{d-1}\phi^{e},\gamma_{1},\cdots,\gamma_{m}\rangle_{g,m+1,\beta}\nonumber\\
&&+ \sum_{\beta_{1}+\beta_{2}=\beta}\langle \gamma, \gamma_{m+1}\psi^{d},\Delta^{a}\rangle_{0,3,\beta_{1}}
\langle \Delta_{a}\phi^{e-1},\gamma_{1},\cdots,\gamma_{m}\rangle_{g,m+1,\beta_{2}}.
\eea
We prove this by induction on $e$. When $e=0$ this reduces to the divisor equation of gravitational descendants.
By (4) of \cite{KM}, we have
\bea\label{17}
&&\langle\gamma,\gamma_{1},\cdots,\gamma_{m},\gamma_{m+1}\psi^{d-1}\phi^{e+1}\rangle_{g,m+2,\beta}\nonumber\\
&=& \langle\gamma,\gamma_{1},\cdots,\gamma_{m},\gamma_{m+1}\psi^{d}\phi^{e}\rangle_{g,m+2,\beta}-
\sum_{\beta_{1}+\beta_{2}=\beta}\langle  \gamma_{m+1}\psi^{d-1},\Delta^{a}\rangle_{0,2,\beta_{1}}
\langle \Delta_{a}\phi^{e},\gamma,\gamma_{1},\cdots,\gamma_{m}\rangle_{g,m+2,\beta_{2}}\nonumber\\
&=& (\gamma \cap \beta)\langle\gamma_{1},\cdots,\gamma_{m},\gamma_{m+1}\psi^{d}\phi^{e}\rangle_{g,m+1,\beta}
+\sum_{\beta_{1}+\beta_{2}=\beta}\langle \gamma, \gamma_{m+1}\psi^{d},\Delta^{a}\rangle_{0,3,\beta_{1}}
\langle \Delta_{a}\phi^{e-1},\gamma_{1},\cdots,\gamma_{m}\rangle_{g,m+1,\beta_{2}}\nonumber\\
&&+ \langle (\gamma\cup\gamma_{m+1})\psi^{d-1}\phi^{e},\gamma_{1},\cdots,\gamma_{m}\rangle_{g,m+1,\beta}
-\sum_{\beta_{1}+\beta_{2}=\beta}\langle  \gamma_{m+1}\psi^{d-1},\Delta^{a}\rangle_{0,2,\beta_{1}}
\langle \Delta_{a}\phi^{e},\gamma,\gamma_{1},\cdots,\gamma_{m}\rangle_{g,m+2,\beta_{2}},\nonumber\\
&&\eea
but
\bea\label{18}
&&\sum_{\beta_{1}+\beta_{2}=\beta}\langle \gamma, \gamma_{m+1}\psi^{d},\Delta^{a}\rangle_{0,3,\beta_{1}}
\langle \Delta_{a}\phi^{e-1},\gamma_{1},\cdots,\gamma_{m}\rangle_{g,m+1,\beta_{2}}\nonumber\\
&&-\sum_{\beta_{1}+\beta_{2}=\beta}\langle  \gamma_{m+1}\psi^{d-1},\Delta^{a}\rangle_{0,2,\beta_{1}}
\langle \Delta_{a}\phi^{e},\gamma,\gamma_{1},\cdots,\gamma_{m}\rangle_{g,m+2,\beta_{2}}\nonumber\\
&=&\sum_{\beta_{1}+\beta_{2}=\beta}\langle \gamma, \gamma_{m+1}\psi^{d},\Delta^{a}\rangle_{0,3,\beta_{1}}
\langle \Delta_{a}\phi^{e-1},\gamma_{1},\cdots,\gamma_{m}\rangle_{g,m+1,\beta_{2}}\nonumber\\
&&-\sum_{\beta_{1}+\beta_{2}=\beta}\langle  \gamma_{m+1}\psi^{d-1},\Delta^{a}\rangle_{0,2,\beta_{1}}(\gamma\cap\beta_{2})
\langle \Delta_{a}\phi^{e},\gamma_{1},\cdots,\gamma_{m}\rangle_{g,m+1,\beta_{2}}\nonumber\\
&&-\sum_{\beta_{1}+\beta_{2}+\beta_{3}=\beta}\langle  \gamma_{m+1}\psi^{d-1},\Delta^{a}\rangle_{0,2,\beta_{1}}
\langle \Delta_{a},\gamma,\Delta^{b}\rangle_{0,3,\beta_{3}}\langle \Delta_{b}\phi^{e-1},\gamma_{1},\cdots,\gamma_{m}\rangle_{g,m+1,\beta_{2}}\nonumber\\
&=&-\sum_{\beta_{1}+\beta_{2}=\beta}\langle  \gamma_{m+1}\psi^{d-1},\Delta^{a}\rangle_{0,2,\beta_{1}}(\gamma\cap\beta_{2})
\langle \Delta_{a}\phi^{e},\gamma_{1},\cdots,\gamma_{m}\rangle_{g,m+1,\beta_{2}},
\eea
in the last equality we have used (4a) in \cite{KM}. Thus
\bea\label{19}
&&\langle\gamma,\gamma_{1},\cdots,\gamma_{m},\gamma_{m+1}\psi^{d-1}\phi^{e+1}\rangle_{g,m+2,\beta}\nonumber\\
&=&(\gamma \cap \beta)\langle\gamma_{1},\cdots,\gamma_{m},\gamma_{m+1}\psi^{d}\phi^{e}\rangle_{g,m+1,\beta}
+ \langle (\gamma\cup\gamma_{m+1})\psi^{d-1}\phi^{e},\gamma_{1},\cdots,\gamma_{m}\rangle_{g,m+1,\beta}\nonumber\\
&&-\sum_{\beta_{1}+\beta_{2}=\beta}\langle  \gamma_{m+1}\psi^{d-1},\Delta^{a}\rangle_{0,2,\beta_{1}}(\gamma\cap\beta_{2})
\langle \Delta_{a}\phi^{e},\gamma_{1},\cdots,\gamma_{m}\rangle_{g,m+1,\beta_{2}}.
\eea
On the other hand, still by (4) of \cite{KM} and the divisor equation for gravitational descendants, for $d\geq 2$ we have
\bea\label{20}
&&(\gamma \cap \beta)\langle \gamma_{1},\cdots,\gamma_{m},\gamma_{m+1}\psi^{d-1}\phi^{e+1}\rangle_{g,m+1,\beta}
+\langle \gamma\cup\gamma_{m+1}\psi^{d-2}\phi^{e+1},\gamma_{1},\cdots,\gamma_{m}\rangle_{g,m+1,\beta}\nonumber\\
&&+ \sum_{\beta_{1}+\beta_{2}=\beta}\langle \gamma, \gamma_{m+1}\psi^{d-1},\Delta^{a}\rangle_{0,3,\beta_{1}}
\langle \Delta_{a}\phi^{e},\gamma_{1},\cdots,\gamma_{m}\rangle_{g,m+1,\beta_{2}}\nonumber\\
&=& (\gamma \cap \beta)\Big(\langle \gamma_{1},\cdots,\gamma_{m},\gamma_{m+1}\psi^{d}\phi^{e}\rangle_{g,m+1,\beta}
-\sum_{\beta_{1}+\beta_{2}=\beta}\langle  \gamma_{m+1}\psi^{d-1},\Delta^{a}\rangle_{0,2,beta_{1}}
\langle \Delta_{a}\phi^{e},\gamma_{1},\cdots,\gamma_{m}\rangle_{g,m+1,\beta_{2}}\Big)\nonumber\\
&&+\langle \gamma\cup\gamma_{m+1}\psi^{d-2}\phi^{e+1},\gamma_{1},\cdots,\gamma_{m}\rangle_{g,m+1,\beta}\nonumber\\
&&+ \sum_{\beta_{1}+\beta_{2}=\beta}(\gamma\cap\beta_{1})\langle \gamma_{m+1}\psi^{d-1},\Delta^{a}\rangle_{0,2,\beta_{1}}
\langle \Delta_{a}\phi^{e},\gamma_{1},\cdots,\gamma_{m}\rangle_{g,m+1,\beta_{2}}\nonumber\\
&&+\sum_{\beta_{1}+\beta_{2}=\beta}\langle (\gamma\cup\gamma_{m+1})\psi^{d-2},\Delta^{a}\rangle_{0,2,\beta_{1}}
\langle \Delta_{a}\phi^{e},\gamma_{1},\cdots,\gamma_{m}\rangle_{g,m+1,\beta_{2}}\nonumber\\
&=&(\gamma \cap \beta)\langle \gamma_{1},\cdots,\gamma_{m},\gamma_{m+1}\psi^{d}\phi^{e}\rangle_{g,m+1,\beta}
-\sum_{\beta_{1}+\beta_{2}=\beta}\langle  \gamma_{m+1}\psi^{d-1},\Delta^{a}\rangle_{0,2,\beta_{1}}(\gamma\cap\beta_{2})
\langle \Delta_{a}\phi^{e},\gamma_{1},\cdots,\gamma_{m}\rangle_{g,m+1,\beta_{2}}\nonumber\\
&&+\Big(\langle \gamma\cup\gamma_{m+1}\psi^{d-2}\phi^{e+1},\gamma_{1},\cdots,\gamma_{m}\rangle_{g,m+1,\beta}
+\sum_{\beta_{1}+\beta_{2}=\beta}\langle (\gamma\cup\gamma_{m+1})\psi^{d-2},\Delta^{a}\rangle_{0,2,\beta_{1}}
\langle \Delta_{a}\phi^{e},\gamma_{1},\cdots,\gamma_{m}\rangle_{g,m+1,\beta_{2}}\Big).\nonumber\\
&&
\eea
Comparing (\ref{19}) and (\ref{20}) we see that (\ref{17}) holds for $e+1$. The case of $d=1$ is similar.\hfill\qedsymbol\\

We need only the ancestor divisor equation. Let $\gamma_{0},\gamma_{1},\cdots, \gamma_{n-1}$ be a basis of $H^{*}(V)$, with $\gamma_{1},\cdots, \gamma_{s}$ a basis of $H^{2}(V)$,
such that $\langle \gamma_{i}, \beta_{j}\rangle=\delta_{i,j}$ for $1\leq i,j\leq s$. Let $q_{1},\cdots, q_{s}$ be the corresponding K\"{a}hler parameter. For $d=d_{1}\beta_{1}+\cdots d_{s}\beta_{s}$, write $q^{d}=q_{1}^{d_{1}}\cdots q_{s}^{d_{s}}$.
Let $\tau=t^{1}\gamma_{1}+\cdots+t^{s}\gamma_{s}$.\\
 For $0\leq i,j\leq n-1$,
\ben
&&\sum_{m\geq 0}\sum_{d\geq 0}\frac{1}{m!}\langle \tau,\cdots,\tau,\gamma_{i},\gamma_{j},\frac{\widetilde{t^{\alpha}}\gamma_{\alpha}}{z-\phi}\rangle_{0,m+3,d}q^{d}\\
&=& \exp \Big(\sum_{k=1}^{s}t^{k}q_{k}\frac{\partial}{\partial q_{k}}
 + \sum_{d\geq 0} \langle \frac{\tilde{t^{\alpha}}\gamma_{\alpha}}{z} ,\tau,\gamma^{\beta}\frac{\partial}{\partial \widetilde{t^{\beta}}}
 \rangle_{0,3,d}q^{d}\Big)\sum_{d\geq 0} \langle \frac{\widetilde{t^{\alpha}}\gamma_{\alpha}}{z} ,\gamma_{i},\gamma_{j}
 \rangle_{0,3,d}q^{d},
\een
where the Greek subscripts suggest the Einstein's convention, i.e. a summation over $\{0,\cdots,n-1\}$.\\
We denote the small quantum product by $\star$, then
\ben
\sum_{m\geq 0}\sum_{d\geq 0}\frac{1}{m!}\langle \tau,\cdots,\tau,\gamma_{i},\gamma_{j},\frac{\gamma_{\alpha}}{z-\phi}\rangle_{0,m+3,d}q^{d}\gamma^{\alpha}=
 \exp \Big(\sum_{k=1}^{s}t^{k}q_{k}\frac{\partial}{\partial q_{k}}
 + \frac{\tau}{z}\star\Big)\sum_{d\geq 0} \frac{\gamma_{i}\star\gamma_{j}}{z}.
\een
When $\gamma_{i}$ and $\gamma_{j}$ are also divisor classes, summing over all of them,  we obtain
\ben
z+\tau+\sum_{m\geq 0}\sum_{d\geq 0}\frac{1}{m!}\langle \tau,\cdots,\tau,\frac{\gamma_{\alpha}}{z-\phi}\rangle_{0,m+1,d}q^{d}\gamma^{\alpha}=
 z\exp \Big(\sum_{k=1}^{s}t^{k}q_{k}\frac{\partial}{\partial q_{k}}
 + \frac{\tau}{z}\star\Big)1.
\een

\begin{remark}\label{28}
For a hard-Lefschetzian orbifold $\mathcal{X}$ and a crepant resolution $Y$, assuming that the crepant resolution conjecture holds for the ancestor correlators (see the introduction for the meaning of the analytic continuations), one can determine the ancestor correlators (in the absolute stable range) of $Y$
from those of $\mathcal{X}$ as follows. \\

We use the same notations as above. For simplicity we assume all the classes $\beta_{1},\cdots,\beta_{s}\in H_{2}(Y)$ are contracted in $\pi:Y\rightarrow X$, where $X$ is the coarse moduli of $\mathcal{X}$. To determine the correlators $\langle \tau_{0,k_{1}}\gamma_{l_1},\cdots, \tau_{0,k_{m}}\gamma_{l_m}\rangle_{g,n,d}$ is equivalent to determine the series
\bea\label{29}
\sum_{d\geq 0}\langle \tau_{0,k_{1}}\gamma_{1},\cdots, \tau_{0,k_{m}}\gamma_{l_m}\rangle_{g,m,d}q^{d},
\eea
 thus is equivalent to determine its Taylor expansion at $(q_{1},\cdots,q_{s})=(\omega_{1},\cdots,\omega_{s})$, where $\omega_{i}$ is the number at which $q_{i}$ takes values (after analytic continuations) to make correspondences with the correlators of $\mathcal{X}$. The value of (\ref{29}) at $(q_{1},\cdots,q_{s})=(\omega_{1},\cdots,\omega_{s})$, where $\omega_{i}$ is determined by the corresponding correlator of $\mathcal{X}$.
Furthermore, by the ancestor divisor equation, we have for $1\leq i\leq s$
\ben
&&\sum_{d\geq 0}\langle\tau_{0,0}\gamma_{i} \tau_{0,k_{1}}\gamma_{l_1},\cdots, \tau_{0,k_{m}}\gamma_{l_m}\rangle_{g,m,d}q^{d}\\
&=&q_{i}\frac{d}{d q_{i}}\sum_{d\geq 0}\langle \tau_{0,k_{1}}\gamma_{1},\cdots, \tau_{0,k_{m}}\gamma_{m}\rangle_{g,m,d}q^{d}\\
&&+\sum_{j=1}^{m}\sum_{d\geq 0}\langle \tau_{0,k_{1}}\gamma_{l_1},\cdots, \tau_{0,k_{j}-1}(\gamma_{i}\star\gamma_{l_j})\cdots \tau_{0,k_{m}}\gamma_{l_m}\rangle_{g,m,d}q^{d}.
\een
By the crepant resolution conjecture, taking $(q_{1},\cdots,q_{s})=(\omega_{1},\cdots,\omega_{s})$, the LHS and the second group of terms of the RHS are determined by the correlators of $\mathcal{X}$. The operator $q_{i}\frac{d}{d q_{i}}=(q_{i}-\omega_{i})\frac{d}{d (q_{i}-\omega_{i})}+\omega_{i}\frac{d}{d (q_{i}-\omega_{i})}$ gives the coefficient of $q_{i}-\omega_{i}$ in the Taylor expansion of (\ref{29}) at $(q_{1},\cdots,q_{s})=(\omega_{1},\cdots,\omega_{s})$. In the same way we can determine the other coefficients inductively.\hfill\qedsymbol
\end{remark}

\subsection{Analytic continuation of the ancestor $J$-function of $\widehat{[\mathbb{C}^{2}/\mathbb{Z}_{2}]}$}
For $Y=\widehat{[\mathbb{C}^{2}/\mathbb{Z}_{2}]}$, let $\gamma_{0}=1$, $\gamma_{1}$ be the equivariant chern class of the equivariant line bundle with weight $-\lambda_{1}$ and $-\lambda_{2}$ at the two fixed points, as in \cite{BGr}. We have seen in the last section that
\bea\label{1}
J_{Y}^{An}/e^{\frac{t^{0}}{z}}=z\exp \Big(t^{1}q\frac{\partial}{\partial q}
 + \frac{t^{1}\gamma_{1}}{z}\star\Big)1,
\eea
where (cf.\cite{BGr})
\ben
\gamma_{1}\star \gamma_{1}=-\lambda_{1}\lambda_{2}\gamma_{0}-\frac{1+q}{1-q}(\lambda_{1}+\lambda_{2})\gamma_{1}.
\een
On the other hand, the $I$-function of $\mathcal{X}=[\mathbb{C}^{2}/\mathbb{Z}_{2}]$ is given by
\bea\label{2}
I_{\mathcal{X}}/e^{\frac{x_0}{z}}&=&z\delta_{0}+x_{1}\delta_{\frac{1}{2}}+\frac{\lambda_1\lambda_2 x_{1}^{2}}{2z}\delta_{0}+\frac{\lambda_1\lambda_2}{z}\sum_{k\geq 2}\prod_{1\leq r\leq k-1}(\frac{\lambda_{1}}{z}-r)(\frac{\lambda_{2}}{z}-r)\frac{x_{1}^{2k}}{(2k)!}\delta_{0}\nonumber\\
&&+\sum_{k\geq 1}\prod_{0\leq r\leq k-1}(\frac{\lambda_{1}}{z}-(r+\frac{1}{2}))(\frac{\lambda_{2}}{z}-(r+\frac{1}{2}))\frac{x_{1}^{2k+1}}{(2k+1)!}\delta_{\frac{1}{2}},
\eea
and after the change of variables
\ben
\tau^{0}&=&x_{0},\\
\tau^{1}&=&\sum_{k\geq 0}\frac{((k-\frac{1}{2})!)^{2}}{(2k+1)!}x_{1}^{2k+1}=2\arcsin (\frac{x_{1}}{2}),
\een
we have the equality $I_{\mathcal{X}}=J_{\mathcal{X}}$, where $J_{\mathcal{X}}$ is given by
\ben
J_{\mathcal{X}}(\tau,z)/e^{\frac{\tau^0}{z}}&=& z+\tau^{1}\delta_{\frac{1}{2}}+\sum_{m\geq 0}\sum_{k=0}^{1}\frac{1}{m!}
\langle \tau^{1}\delta_{\frac{1}{2}},\cdots,\tau^{1}\delta_{\frac{1}{2}},\frac{\delta_{k}}{z-\psi}\rangle_{0,m+1}^{\mathcal{X}}\delta^{k}\\
&=& \Big[z+\sum_{r\geq 1}\frac{2\lambda_1\lambda_2}{z^{r}}\sum_{m\geq 2}\frac{(\tau^{1})^{m}}{m!}
\langle \delta_{\frac{1}{2}},\cdots,\delta_{\frac{1}{2}},\delta_{0}\psi^{r-1}\rangle_{0,m+1}^{\mathcal{X}}\Big]\delta_{0}\\
&&+\Big[\tau^{1}+\sum_{r\geq 1}\frac{2}{z^{r}}\sum_{m\geq 2}\frac{(\tau^{1})^{m}}{m!}
\langle \delta_{\frac{1}{2}},\cdots,\delta_{\frac{1}{2}},\delta_{\frac{1}{2}}\psi^{r-1}\rangle_{0,m+1}^{\mathcal{X}}
\Big]\delta_{\frac{1}{2}}
\een
since $\delta^{0}=2\lambda_1\lambda_2\delta_{0}$, $\delta^{\frac{1}{2}}=2\delta_{\frac{1}{2}}$. Our goal is to show
\begin{theorem}\label{26}
 $J_{Y}^{An}=J_{\mathcal{X}}$ after the continuation of $q$ from 0 to $-1$ (along the negative real axis), and the change of variables and cohomology classes $t^{0}=\tau^{0}$, $t^{1}=\sqrt{-1}\tau$, $\gamma_{0}=\delta_{0}=1$,$\gamma_{1}=-\sqrt{-1}\delta_{\frac{1}{2}}$.
\end{theorem}
Thus it suffices to compare the RHS of (\ref{1}) after analytic continuations and the RHS of (\ref{2}). For example, we expand the RHS of (\ref{1}) up to $(t^{1})^{4}$,
\ben
J_{Y}^{An}/e^{\frac{t^{0}}{z}}&=&\Big[-\frac{\lambda_{1}\lambda_{2}(t^{1})^{2}}{z^{2}}
+\frac{1+q}{1-q}(\lambda_{1}+\lambda_{2})\lambda_{1}\lambda_{2}\frac{(t^{1})^{3}}{z^{3}}\\
&&+(t^{1})^{4}\Big(
 \frac{4q}{(1-q)^{2}}\frac{(\lambda_{1}+\lambda_{2})\lambda_{1}\lambda_{2}}{z^{3}}-\frac{\big((\frac{1+q}{1-q})^{2}
 (\lambda_{1}+\lambda_{2})^{2}-\lambda_{1}\lambda_{2}\big)\lambda_{1}\lambda_{2}}{z^{4}}\Big)
\Big]\gamma_{0}\\
&&+\Big[\frac{t^{1}}{z} -\frac{1+q}{1-q}\frac{(\lambda_{1}+\lambda_{2})(t^{1})^{2}}{z^{2}}
+\Big(-\frac{2q}{(1-q)^{2}}\frac{(\lambda_{1}+\lambda_{2})(t^{1})^{3}}{z^{2}} +\big((\frac{1+q}{1-q})^{2}(\lambda_{1}+\lambda_{2})^{2}-\lambda_{1}\lambda_{2}\big)\frac{(t^{1})^{3}}{z^{3}}
 \Big)\\
&&+(t^{1})^{4}\Big(-\frac{2q(1+q)}{(1-q)^{3}}\frac{\lambda_{1}+\lambda_{2}}{z^{2}}
 +\frac{6q(1+q)}{(1-q)^{3}}\frac{(\lambda_{1}+\lambda_{2})^{2}}{z^{3}}\\
 &&+\frac{\frac{2(1+q)}{1-q}(\lambda_{1}+\lambda_{2})\lambda_{1}\lambda_{2}-(\frac{1+q}{1-q})^{3}(\lambda_{1}+\lambda_{2})^{3}}{z^{4}}
 \Big)\Big]\gamma_{1}+O(t^{5}).
\een
A simple and crucial observation:\\

 \emph{In this expression, analytic continuation from $q=0$  to $q=-1$ along the negative real axis means no other than directly take $q=-1$.} \hfill ($*$)\\

Thus it's not hard to see that, after this specification of $q$ and the change of variables
\bea\label{20}
t^{1}=2\sqrt{-1}\arcsin (\frac{x_{1}}{2}),
\eea
$J_{Y}^{An}/e^{\frac{t^{0}}{z}}$ coincides with the RHS of (\ref{2}) up to $x_{1}^{4}$. In general, we cannot give a closed formula for the RHS of
(\ref{1}). However, setting $J_{Y}^{An}/e^{\frac{t^{0}}{z}}=\Phi_{0}\gamma_{0}+\Phi_{1}\gamma_{1}$, we have a \emph{quantum differential equation} (QDE for short)
\bea\label{3}
q\frac{\partial}{\partial q}\left(\begin{array}{c}
\Phi_{0}\\
\Phi_{1}\end{array}
\right)=\left(\begin{array}{cc}
 \frac{\partial}{\partial t^{1}} & \frac{\lambda_{1}\lambda_{2}}{z}\\
 -\frac{1}{z}  & \frac{\partial}{\partial t^{1}} + \frac{(1+q)(\lambda_{1}+\lambda_{2})}{(1-q)z}
\end{array}
\right)\left(\begin{array}{c}
\Phi_{0}\\
\Phi_{1}\end{array}
\right).
\eea
Along every characteristic curve $q=\mu e^{-t^{1}}$ of the differential operator $q\frac{\partial}{\partial q}-\frac{\partial}{\partial t^{1}}$, having  known the initial value at $t^{1}=0$
\bea\label{21}
J_{Y}^{An}/e^{\frac{t^{0}}{z}}|_{t^{1}=0}=z,
\eea
 we can integrate the corresponding ODE to obtain the value of $J_{Y}^{An}/e^{\frac{t^{0}}{z}}$ at $q=-1$. Varying $\mu$, by the observation ($*$), we obtain the analytic continuation we need. This is our strategy.
\subsection{Integrating the QDE}
From (\ref{3}) we easily get
\ben
q^{2}\frac{d^{2}}{d q^{2}}\Phi_{0}+\Big(1-\frac{1+q}{1-q}\frac{\lambda_{1}+\lambda_{2}}{z}\Big)q\frac{d}{d q}\Phi_{0}+\frac{\lambda_{1}\lambda_{2}}{z^{2}}\Phi_{0}=0.
\een

Make the change of variables (the \emph{mirror map})

\ben
q=-e^{i\cdot2\arcsin \frac{y}{2}}=\frac{y^{2}}{2}-1-y\sqrt{\frac{y^{2}}{4}-1},
\een
we obtain
\ben
(\frac{y^{2}}{4}-1)\frac{d^{2}}{dy^{2}}\Phi_{0}+(\frac{1}{4}+\frac{\lambda_{1}+\lambda_{2}}{2z})y\frac{d}{dy}\Phi_{0}
+\frac{\lambda_{1}\lambda_{2}}{z^{2}}\Phi_{0}=0,
\een
i.e.,
\ben
\Big[(\frac{y}{2}\frac{d}{dy}+\frac{\lambda_{1}}{z})(\frac{y}{2}\frac{d}{dy}+\frac{\lambda_{2}}{z})-\frac{d^{2}}{dy^{2}}\Big]\Phi_{0}=0.
\een
Note that $y=0$ is an ordinary point of this differential equation. Let
\ben
\Phi_{0}=\sum_{k\geq 0}c_{k}y^{k},
\een
we have
\ben
\frac{k(k-1)}{4}c_{k}-(k+2)(k+1)c_{k+2}+(\frac{1}{4}+\frac{\lambda_{1}+\lambda_{2}}{2z})kc_{k}+\frac{\lambda_{1}\lambda_{2}}{z^{2}}c_{k}=0,
\een
which implies
\ben
c_{k+2}=\frac{(\frac{k}{2}+\frac{\lambda_{1}}{z})(\frac{k}{2}+\frac{\lambda_{2}}{z})}{(k+2)(k+1)}c_{k}.
\een
Thus
\ben
\Phi_{0}=c_{0}\sum_{k\geq 0}\frac{y^{2k}}{(2k)!}\prod_{r=1}^{k}(\frac{\lambda_{1}}{z}+r-1)(\frac{\lambda_{2}}{z}+r-1)
+c_{1}\sum_{k\geq 0}\frac{y^{2k+1}}{(2k+1)!}\prod_{r=1}^{k}(\frac{\lambda_{1}}{z}+r-\frac{1}{2})(\frac{\lambda_{2}}{z}+r-\frac{1}{2}),
\een
and
\ben
\Phi_{1}&=&-\frac{z}{\lambda_{1}\lambda_{2}}\sqrt{\frac{y^{2}}{4}-1}\Big[
c_{0}\sum_{k\geq 1}\frac{y^{2k-1}}{(2k-1)!}\prod_{r=1}^{k}(\frac{\lambda_{1}}{z}+r-1)(\frac{\lambda_{2}}{z}+r-1)\\
&&+c_{1}\sum_{k\geq 0}\frac{y^{2k}}{(2k)!}\prod_{r=1}^{k}(\frac{\lambda_{1}}{z}+r-\frac{1}{2})(\frac{\lambda_{2}}{z}+r-\frac{1}{2})\Big].
\een
Since $\Phi_{1}|_{y=x}=0$, we have
\ben
c_{0}&=&C\cdot \sum_{k\geq 0}\frac{x^{2k}}{(2k)!}\prod_{r=1}^{k}(\frac{\lambda_{1}}{z}+r-\frac{1}{2})(\frac{\lambda_{2}}{z}+r-\frac{1}{2}),\\
c_{1}&=&-C\cdot \sum_{k\geq 1}\frac{x^{2k-1}}{(2k-1)!}\prod_{r=1}^{k}(\frac{\lambda_{1}}{z}+r-1)(\frac{\lambda_{2}}{z}+r-1),
\een
where $C=C(x;\lambda_{1},\lambda_{2},z)$ is a function to be determined.
Since $\Phi_{0}|_{y=x}=z$, we have
\ben
C&=& z\Big/ \Big[\Big(\sum_{k\geq 0}\frac{x^{2k}}{(2k)!}\prod_{r=1}^{k}(\frac{\lambda_{1}}{z}+r-\frac{1}{2})(\frac{\lambda_{2}}{z}+r-\frac{1}{2})\Big)
\Big(\sum_{k\geq 0}\frac{x^{2k}}{(2k)!}\prod_{r=1}^{k}(\frac{\lambda_{1}}{z}+r-1)(\frac{\lambda_{2}}{z}+r-1)\Big)\\
&&-\Big(\sum_{k\geq 1}\frac{x^{2k-1}}{(2k-1)!}\prod_{r=1}^{k}(\frac{\lambda_{1}}{z}+r-1)(\frac{\lambda_{2}}{z}+r-1)\Big)
\Big(\sum_{k\geq 0}\frac{x^{2k+1}}{(2k+1)!}\prod_{r=1}^{k}(\frac{\lambda_{1}}{z}+r-\frac{1}{2})(\frac{\lambda_{2}}{z}+r-\frac{1}{2})\Big)
\Big].
\een
Consequently,
\ben
\Phi_{0}|_{y=0}&=&C\cdot \sum_{k\geq 0}\frac{x^{2k}}{(2k)!}\prod_{r=1}^{k}(\frac{\lambda_{1}}{z}+r-\frac{1}{2})(\frac{\lambda_{2}}{z}+r-\frac{1}{2}),\\
\Phi_{1}|_{y=0}&=&\frac{z\sqrt{-1}}{\lambda_{1}\lambda_{2}}C\cdot \sum_{k\geq 1}\frac{x^{2k-1}}{(2k-1)!}\prod_{r=1}^{k}(\frac{\lambda_{1}}{z}+r-1)(\frac{\lambda_{2}}{z}+r-1).
\een

We need the following
\begin{lemma}
\bea\label{22}
&&C\cdot \sum_{k\geq 0}\frac{x^{2k}}{(2k)!}\prod_{r=1}^{k}(\frac{\lambda_{1}}{z}+r-\frac{1}{2})(\frac{\lambda_{2}}{z}+r-\frac{1}{2})\nonumber\\
&=& z \sum_{k\geq 0}\frac{x^{2k}}{(2k)!}\prod_{r=1}^{k}(\frac{\lambda_{1}}{z}-r+1)(\frac{\lambda_{2}}{z}-r+1),
\eea
and
\bea\label{23}
&&\frac{z}{\lambda_{1}\lambda_{2}}C\cdot \sum_{k\geq 1}\frac{x^{2k-1}}{(2k-1)!}\prod_{r=1}^{k}(\frac{\lambda_{1}}{z}+r-1)(\frac{\lambda_{2}}{z}+r-1)\nonumber\\
&=&  \sum_{k\geq 0}\frac{x^{2k+1}}{(2k+1)!}\prod_{r=1}^{k}(\frac{\lambda_{1}}{z}-r+\frac{1}{2})(\frac{\lambda_{2}}{z}-r+\frac{1}{2}).
\eea
\end{lemma}
\emph{Proof}: Let
\ben
f_{1}(a,b,x)&=&\sum_{k\geq 0}\frac{x^{2k}}{(2k)!}\prod_{r=1}^{k}(a+r-1)(b+r-1),\\
f_{2}(a,b,x)&=&\sum_{k\geq 0}\frac{x^{2k+1}}{(2k+1)!}\prod_{r=1}^{k}(a+r-\frac{1}{2})(b+r-\frac{1}{2}),\\
f_{3}(a,b,x)&=&\sum_{k\geq 0}\frac{x^{2k}}{(2k)!}\prod_{r=1}^{k}(a-r+1)(b-r+1),\\
f_{4}(a,b,x)&=&\sum_{k\geq 0}\frac{x^{2k+1}}{(2k+1)!}\prod_{r=1}^{k}(a-r+\frac{1}{2})(b-r+\frac{1}{2}).
\een
Then
\ben
C=\frac{z}{f_{2}^{\prime}(\frac{\lambda_{1}}{z},\frac{\lambda_{2}}{z},x)f_{1}(\frac{\lambda_{1}}{z},\frac{\lambda_{2}}{z},x)
-f_{1}^{\prime}(\frac{\lambda_{1}}{z},\frac{\lambda_{2}}{z},x)f_{2}(\frac{\lambda_{1}}{z},\frac{\lambda_{2}}{z},x)},
\een
where $\prime=\frac{\partial}{\partial x}$. It reduces to show
\bea\label{24}
(f_{2}^{\prime}f_{1}-f_{1}^{\prime}f_{2})f_{3}=f_{2}^{\prime},
\eea
and
\bea\label{25}
ab(f_{2}^{\prime}f_{1}-f_{1}^{\prime}f_{2})f_{4}=f_{1}^{\prime}.
\eea
From the equation
\ben
(\frac{x^{2}}{4}-1)\frac{d^{2}}{dx^{2}}f_{2}+(\frac{1}{4}+\frac{a+b}{2})x\frac{d}{dx}f_{2}
+ab f_{2}=0,
\een
we easily obtain
\ben
\Big[(\frac{x^{2}}{4}-1)\frac{d^{2}}{dx^{2}}+(\frac{3}{4}+\frac{a+b}{2})x\frac{d}{dx}
+(a+\frac{1}{2})(b+\frac{1}{2})\Big]f_{2}^{\prime}=0.
\een
 We need only to show that the LHS of (\ref{24}) satisfies the same differential equation, for it is easily seen that the two handsides of (\ref{24}) have the same initial values. So we come to show
\ben
\Big[(\frac{x^{2}}{4}-1)\frac{d^{2}}{dx^{2}}+(\frac{3}{4}+\frac{a+b}{2})x\frac{d}{dx}
+(a+\frac{1}{2})(b+\frac{1}{2})\Big](f_{2}^{\prime}f_{1}-f_{1}^{\prime}f_{2})f_{3}=0.
\een
But
\ben
&&\Big[(\frac{x^{2}}{4}-1)\frac{d^{2}}{dx^{2}}+(\frac{3}{4}+\frac{a+b}{2})x\frac{d}{dx}
+(a+\frac{1}{2})(b+\frac{1}{2})\Big](f_{2}^{\prime}f_{1}-f_{1}^{\prime}f_{2})f_{3}\\
&=& (\frac{x^{2}}{4}-1)\Big[(f_{2}^{\prime\prime\prime}f_{1}+f_{2}^{\prime\prime}f_{1}^{\prime}-
f_{2}^{\prime}f_{1}^{\prime\prime}-f_{2}f_{1}^{\prime\prime\prime})f_{3}\\
&&+2(f_{2}^{\prime\prime}f_{1}-f_{2}f_{1}^{\prime\prime})f_{3}^{\prime}+
(f_{2}^{\prime}f_{1}-f_{2}f_{1}^{\prime})f_{3}^{\prime\prime}
\Big]\\
&&+(\frac{3}{4}+\frac{a+b}{2})x\Big[(f_{2}^{\prime\prime}f_{1}-f_{2}f_{1}^{\prime\prime})f_{3}+
(f_{2}^{\prime}f_{1}-f_{2}f_{1}^{\prime})f_{3}^{\prime}
\Big]\\
&&+(a+\frac{1}{2})(b+\frac{1}{2})\Big](f_{2}^{\prime}f_{1}-f_{1}^{\prime}f_{2})f_{3}\\
&=& f_{1}f_{3}\Big[(\frac{x^{2}}{4}-1)f_{2}^{\prime\prime\prime}
+(\frac{3}{4}+\frac{a+b}{2})x(f_{2}^{\prime\prime}+(a+\frac{1}{2})(b+\frac{1}{2})f_{2}^{\prime}
\Big]\\
&&- f_{2}f_{3}\Big[(\frac{x^{2}}{4}-1)f_{1}^{\prime\prime\prime}
+(\frac{3}{4}+\frac{a+b}{2})x(f_{1}^{\prime\prime}+(a+\frac{1}{2})(b+\frac{1}{2})f_{1}^{\prime}
\Big]\\
&&+(\frac{x^{2}}{4}-1)\Big[(f_{2}^{\prime\prime}f_{1}^{\prime}-
f_{2}^{\prime}f_{1}^{\prime\prime})f_{3}+2(f_{2}^{\prime\prime}f_{1}-f_{2}f_{1}^{\prime\prime})f_{3}^{\prime}+
(f_{2}^{\prime}f_{1}-f_{2}f_{1}^{\prime})f_{3}^{\prime\prime}
\Big]\\
&&+(\frac{3}{4}+\frac{a+b}{2})x
(f_{2}^{\prime}f_{1}-f_{2}f_{1}^{\prime})f_{3}^{\prime}\\
&=&(\frac{x^{2}}{4}-1)\Big[(f_{2}^{\prime\prime}f_{1}^{\prime}-
f_{2}^{\prime}f_{1}^{\prime\prime})f_{3}+2(f_{2}^{\prime\prime}f_{1}-f_{2}f_{1}^{\prime\prime})f_{3}^{\prime}+
(f_{2}^{\prime}f_{1}-f_{2}f_{1}^{\prime})f_{3}^{\prime\prime}
\Big]\\
&&+(\frac{3}{4}+\frac{a+b}{2})x
(f_{2}^{\prime}f_{1}-f_{2}f_{1}^{\prime})f_{3}^{\prime}.\\
\een
Note that
\ben
(\frac{x^{2}}{4}-1)f_{2}^{\prime\prime}&=&-(\frac{1}{4}+\frac{a+b}{2})x\frac{d}{dx}f_{2}-ab f_{2},\\
(\frac{x^{2}}{4}-1)f_{1}^{\prime\prime}&=&-(\frac{1}{4}+\frac{a+b}{2})x\frac{d}{dx}f_{1}-ab f_{1},
\een
so we have
\ben
f_{2}^{\prime\prime}f_{1}^{\prime}-f_{2}^{\prime}f_{1}^{\prime\prime}=ab (f_{2}^{\prime}f_{1}-
f_{2}f_{1}^{\prime}),
\een
and
\ben
f_{2}^{\prime\prime}f_{1}-f_{2}f_{1}^{\prime\prime}=-(\frac{1}{4}+\frac{a+b}{2})x(f_{2}^{\prime}f_{1}-
f_{2}f_{1}^{\prime}).
\een
Thus
\ben
&&\Big[(\frac{x^{2}}{4}-1)\frac{d^{2}}{dx^{2}}+(\frac{3}{4}+\frac{a+b}{2})x\frac{d}{dx}
+(a+\frac{1}{2})(b+\frac{1}{2})\Big](f_{2}^{\prime}f_{1}-f_{1}^{\prime}f_{2})f_{3}\\
&=& ab (f_{2}^{\prime}f_{1}-
f_{2}f_{1}^{\prime})f_{3}-2(\frac{1}{4}+\frac{a+b}{2})x(f_{2}^{\prime}f_{1}-
f_{2}f_{1}^{\prime})f_{3}^{\prime}\\
&&+(\frac{x^{2}}{4}-1)(f_{2}^{\prime}f_{1}-
f_{2}f_{1}^{\prime})f_{3}^{\prime\prime}+(\frac{3}{4}+\frac{a+b}{2})x
(f_{2}^{\prime}f_{1}-f_{2}f_{1}^{\prime})f_{3}^{\prime}\\
&=& (f_{2}^{\prime}f_{1}-f_{2}f_{1}^{\prime})\Big(
(\frac{x^{2}}{4}-1)f_{3}^{\prime\prime}+(\frac{1}{4}-\frac{a+b}{2})x f_{3}^{\prime}+abf_{3}.
\Big).
\een
But
\ben
&&x^{2}\Big((\frac{x^{2}}{4}-1)\frac{d^{2}}{dx^{2}}+(\frac{1}{4}-\frac{a+b}{2})x\frac{d}{dx}+ab\Big)\\
&=& x\frac{d}{dx}(x\frac{d}{dx}-1)-x^{2}(a-\frac{x}{2}\frac{d}{dx})(b-\frac{x}{2}\frac{d}{dx}),
\een
which is the hypergeometric operator that annihilate $f_{3}$. Thus we have proven (\ref{24}). Exactly the same procedure applies for (\ref{25}).\hfill\qedsymbol\\

By the above lemma and the discussions in the last section, we complete the proof of theorem \ref{26}. \hfill\qedsymbol\\

\end{appendices}

\hspace{1cm}\footnotesize{Department of Mathematical Sciences, Tsinghua University, Beijing, 100084, China }\\

\hspace{1cm}\footnotesize{\emph{E-mail address}: huxw08@mails.tsinghua.edu.cn}

\end{document}